\documentclass[12pt,a4paper,twoside]{article}

\usepackage[active]{srcltx}
\usepackage[utf8]{inputenc}
\usepackage{hyperref}
\usepackage[english]{babel}
\usepackage{authblk}

\usepackage[a4paper,hcentering,vcentering]{geometry}
\usepackage{url,epsfig,amssymb,amsthm,amsmath}

\numberwithin{equation}{section}

\numberwithin{theorems}{section}

\numberwithin{corollarys}{section}

\newtheorem{definition}{Definition}
\numberwithin{definition}{section}

\newtheorem{lemma}{Lemma}
\numberwithin{lemma}{section}
\newtheorem{proposition}{Proposition}
\numberwithin{proposition}{section}
\newtheorem{theoremx}{Theorem}

\newtheorem{corollaryx}{Corollary}

\newtheorem*{definitionnull}{Definition}

\newcommand{\C}{\ensuremath{\mathbb{C}}}

\newcommand{\N}{\ensuremath{\mathbb{N}}}

\newcommand{\R}{\ensuremath{\mathbb{R}}}
\newcommand{\T}{\ensuremath{\mathbb{T}}}
\newcommand{\Z}{\ensuremath{\mathbb{Z}}}

\begin{document}

\title{\LARGE{\textbf{Asymptotically quasiperiodic solutions for time-dependent Hamiltonians}}}%
\author{Donato Scarcella}
\affil{Université Paris-Dauphine - Ceremade UMR 7534 Place du Maréchal De Lattre De Tassigny, 75016 PARIS.}
\date{\today}%
\maketitle

\begin{abstract}
In 2015, M. Canadell and R. de la Llave consider a time-dependent perturbation of a vector field having an invariant torus supporting quasiperiodic solutions. Under a smallness assumption on the perturbation and assuming the perturbation decays (when $t \to +\infty$) exponentially fast in time, they proved the existence of motions converging in time (when $t \to +\infty$) to quasiperiodic solutions associated with the unperturbed system (asymptotically quasiperiodic solutions). 

In this paper, we generalize this result in the particular case of time-dependent Hamiltonian systems. The exponential decay in time is relaxed (due to the geometrical properties of Hamiltonian systems) and the smallness assumption on the perturbation is removed. 
\end{abstract}

\section{Introduction}

The Kolmogorov-Arnold-Moser (KAM) theory~\cite{{Kol54}, {Arn63a}, {Mos62}} shows the persistence of quasiperiodic solutions in nearly integrable Hamiltonian systems. This theory is interesting, especially for the applications in classical problems in Celestial Mechanics, such as the $n$-body problem. This work analyses non-autonomous perturbations of Hamiltonians having an invariant torus supporting quasiperiodic solutions. We assume that the perturbation satisfies good decay properties when $t \to + \infty$. In this case, we are not looking for quasiperiodic solutions but different types of orbits converging in time (when $t \to + \infty$) to quasiperiodic solutions.

To be more precise, let us introduce the definition of analytic asymptotic KAM torus. Let $B \subset \R^n$ be a ball centred at the origin and $\mathcal{P}$ equal to $\T^n \times B$ or $\T^n$. Moreover, for a given $\upsilon \ge 0$, we introduce the following interval $J_\upsilon = [\upsilon , +\infty) \subset \R$. We consider time-dependent real analytic vector fields $X^t$ and $X^t_0$ on $\mathcal{P}$, for all $t \in J_\upsilon$, and a real analytic embedding $\varphi_0$ from $\T^n$ to $\mathcal{P}$ such that
\begin{align}
\label{hyp1introAKAM}
& \displaystyle \lim_{t \to +\infty}  |X^t - X^t_0|_s = 0,\\
\label{hyp2introAKAM}
& X_0 (\varphi_0(q),t) =\partial_q \varphi_0(q)\omega \hspace{2mm} \mbox{for all $(q, t) \in \T^n \times J_\upsilon$,}
\end{align}
where $\omega \in \R^n$ and $|\cdot |_s$ is the analytic norm (see Appendix \ref{B}). For the sake of clarity, we specify that $\partial_q \varphi_0(q)\omega$ is the element of $\R^{2n}$ having $j$ component equal to 
\begin{equation*}
\Big(\partial_q \varphi_0(q)\omega\Big)_j = \partial_q \varphi_{0,j}(q)\cdot \omega,
\end{equation*}
for all $j = 1, ...,2n$ where $\varphi_0 = (\varphi_{0,1}, ..., \varphi_{0,2n})$. In words, we are considering a time-dependent vector field $X^t$ converging in time to a vector field $X_0^t$ having an invariant torus supporting quasiperiodic dynamics with frequency vector $\omega$.  

\begin{definition}
\label{analasymKAMtorusI}
We assume that $(X, X_0, \varphi_0)$ satisfy~\eqref{hyp1introAKAM} and~\eqref{hyp2introAKAM}. A family of real analytic embeddings $\varphi^t: \T^n \to \mathcal{P}$ is an analytic asymptotic KAM torus associated to $(X, X_0, \varphi_0)$ if there exist $0 < s' \le s$ and $\upsilon' \ge \upsilon \ge 0$ such that
\begin{align*}
&  \lim_{t \to +\infty}  |\varphi^t - \varphi_0|_{s'} = 0,\\
&  X (\varphi(q,t), t) = \partial_q \varphi(q, t) \omega + \partial_t \varphi(q, t), 
\end{align*}
for all $(q, t) \in \T^n \times J_{\upsilon'}$. When $\mathcal{P}$ is a symplectic manifold with $\mathrm{dim} \mathcal{P} = 2n$, then we say that $\varphi^t$ is Lagrangian if $\varphi^t(\T^n)$ is Lagrangian for all $t \in J_{\upsilon'}$. 
\end{definition}

Roughly speaking, an analytic asymptotic KAM torus is a family of embeddings $\varphi^t$ converging in time to the invariant torus $\varphi_0$ associated with $X_0$. Moreover, the dynamics on this family of embeddings converge to the quasiperiodic solutions associated with $X_0$ on $\varphi_0$.
The previous definition is due to M. Canadell and R. de la Llave (see~\cite{CdlL15}). In their paper, they use the expression non-autonomous KAM torus. We prefer analytic asymptotic KAM torus to point out the asymptotic properties of this family of embeddings.  Hence, the following definition is quite natural
\begin{definition}
\label{defqpsolI}
We assume that $(X, X_0, \varphi_0)$ satisfy~\eqref{hyp1introAKAM} and~\eqref{hyp2introAKAM}. An integral curve $g(t)$ of $X$ is an asymptotically quasiperiodic solution associated to $(X, X_0, \varphi_0)$ if  there exists $q \in \T^n$ in such a way that 
\begin{equation*}
\lim_{t \to +\infty} |g(t) - \varphi_0 (q + \omega (t -t_0))|=0. 
\end{equation*}
\end{definition}
It is straightforward to see that if there exists an analytic asymptotic KAM torus associated to $(X, X_0, \varphi_0)$, then we have the existence of asymptotically quasiperiodic solutions associated to $(X, X_0, \varphi_0)$.
We refer to Section \ref{AKAMtori} for a more detailed explanation and a  series of remarks about the previous definitions always due to M. Canadell and R. de la Llave.

The first result about the existence of an analytic asymptotic KAM torus associated with a suitable time-dependent Hamiltonian is due to  A. Fortunati and S. Wiggins~\cite{FW14}.  They consider a real analytic time-dependent Hamiltonian $H^t$ as the sum between an integrable Hamiltonian $h$ plus a time-dependent perturbation $f^t$. They assume $h$ to be non-degenerate and the perturbation decays exponentially fast in time ($|f^t|_s \to_{t \to +\infty} 0$). Suppose the perturbation is sufficiently small. Then, for each invariant torus associated with $h$ supporting quasiperiodic dynamics with a diophantine frequency vector, they prove the existence of an analytic asymptotic KAM torus associated to $(X_H, X_h, \varphi_0)$. For clarity, $X_H$ and $X_h$ are the Hamiltonian systems associated with the Hamiltonians $H$ and $h$. 

About one year later, M. Canadell and R. De la Llave~\cite{CdlL15} publish a result which generalises the one of Fortunati-Wiggins. In this paper, they work with finitely differentiable time-dependent vector fields.  For this reason, let us introduce the definition of $C^\sigma$-asymptotic KAM torus, which is the finitely differentiable version of Definition \ref{analasymKAMtorusI}. 

Given $\sigma \ge 0$ and a positive integer $k \ge 0$, we consider time-dependent vector fields $X^t$ and $X^t_0$ of class $C^{\sigma + k}$ on $\mathcal{P}$, for all $t \in J_{\upsilon}$, and an embedding $\varphi_0$ from $\T^n$ to $\mathcal{P}$ of class $C^\sigma$ such that
\begin{align}
\label{hyp1introKAM}
& \displaystyle \lim_{t \to +\infty} |X^t - X^t_0|_{C^{\sigma + k}} = 0,\\
\label{hyp2introKAM}
& X_0 (\varphi_0(q), t) =\partial_q \varphi_0(q)\omega \hspace{2mm} \mbox{for all $(q, t) \in \T^n \times J_\upsilon$,}
\end{align}
where $\omega \in \R^n$, $C^\sigma$ is the space of the Hölder functions and $|\cdot|_{C^\sigma}$ is the Hölder norm (see Appendix \ref{A}). 
\begin{definition}
\label{asymKAMtoriHol}
We assume that $(X, X_0, \varphi_0)$ satisfy~\eqref{hyp1introKAM} and~\eqref{hyp2introKAM}. A family of $C^\sigma$ embeddings $\varphi^t: \T^n \to \mathcal{P}$ is a $C^\sigma$-asymptotic KAM torus associated to $(X, X_0, \varphi_0)$ if there exists $\upsilon' \ge \upsilon \ge 0$ such that  
\begin{align*}
&  \lim_{t \to +\infty}   |\varphi^t - \varphi_0|_{C^\sigma} = 0,\\
&  X (\varphi(q,t), t) = \partial_q \varphi(q, t) \omega + \partial_t \varphi(q, t), 
\end{align*}
for all $(q, t) \in \T^n \times J_{\upsilon'}$. When $\mathcal{P}$ is a symplectic manifold with $\mathrm{dim} \mathcal{P} = 2n$, then we say that $\varphi^t$ is Lagrangian if $\varphi^t(\T^n)$ is Lagrangian for all $t$.  
\end{definition}

In the work of Canadell-de la Llave, $\mathcal{P}$ is a smooth manifold. They consider a time-dependent vector field $X^t$ converging exponentially fast in time to an autonomous vector field $X_0$ having an invariant torus $\varphi_0$ supporting quasiperiodic solutions of frequency vector $\omega \in \R^n$. They do not assume any non-degeneracy hypothesis on $X_0$ or arithmetic conditions on the frequency vector $\omega$. On the other hand, they need to assume a certain control on the normal dynamics on the invariant torus $\varphi_0$. Then, they prove the existence of a $C^\sigma$-asymptotic KAM torus associated to $(X, X_0, \varphi_0)$. 

In this paper, we generalize the previous result in the particular case of finitely differentiable (or real analytic) Hamiltonian systems. We consider a time-dependent Hamiltonian as a time-dependent perturbation of a Hamiltonian having a Lagrangian invariant torus supporting quasiperiodic solutions. The exponential decay in time is replaced by a more general condition (verified for a suitable polynomial decay in time) and the smallness assumption on the perturbation is removed. Under these hypotheses, we prove the existence of a $C^\sigma$-asymptotic KAM torus (or an analytic asymptotic KAM torus). 

The advantage of working with Hamiltonian systems is that the normal dynamics to the Lagrangian invariant torus associated with the unperturbed system are easy to control. This allows us to relax the hypothesis of exponential decay in time. 
Concerning the smallness condition over the perturbation, we approach the problem from another point of view. We study our system for all time $t$ large enough so that the perturbative terms are sufficiently small and we find the existence of an analytic (or  $C^\sigma$) asymptotic KAM torus defined for all $t$ large. Therefore, we will see that if we prove the existence of an analytic (or  $C^\sigma$) asymptotic KAM torus defined for all $t$ large, then we can extend the set of definition for all $t$ (see Section \ref{AKAMtori}).
The following section contains the main results of this work. But, first, just a few remarks. 

At this point, it seems reasonable to wonder when we have the existence of biasymptotically quasiperiodic solutions. That is, orbits converging to quasiperiodic dynamics of a specific frequency vector $\omega_+ \in \R^n$ in the future $(t \to +\infty)$ and to quasiperiodic dynamics of frequency vector $\omega_- \in \R^n$ in the past $(t \to -\infty)$. It is the subject of another work~\cite{Sca22c}.

The interest in these kinds of perturbations is not artificial. These types of systems would be interesting in astronomy. In a further paper~\cite{Sca22b}, we study the example of a planetary system perturbed by a comet coming from and going back to infinity, asymptotically along a hyperbolic Keplerian orbit.

\section{Results}

This section is divided into two subsections. This is because we prove the same results in the finitely differentiable and analytical case. We look at time-dependent perturbations of suitable time-dependent Hamiltonians and constant vector fields on the torus. We recall that $B \subset \R^n$ is a ball around the origin and, for given $\upsilon \ge 0$, we have the following interval $J_\upsilon = [\upsilon , +\infty) \subset \R$.
In this work, we consider time-dependent Hamiltonians of the form
\begin{equation*}
H:\T^n \times B \times J_0 \to \R, \quad H(q,p,t) = h(q,p,t) + f(q,p,t),
\end{equation*}
where $h$ is in the $\omega$-Kolmogorov normal form, with $\omega\in\R^n$. This means that, for all $(q,t) \in \T^n \times J_0$
\begin{eqnarray*}
h(q,0,t) = c, \quad \partial_p h(q,0,t) = \omega,
\end{eqnarray*}
for some $c \in \R$. Letting $\varphi_0$ be the following trivial embedding $\varphi_0:\T^n \to \T^n \times B$ with $\varphi_0(q) = (q,0)$,
then $\varphi_0$ is a Lagrangian invariant torus for $X_h$ supporting quasiperiodic dynamics with frequency vector $\omega$. We define $\mathcal{K}_\omega$ as the set of the Hamiltonians $h:\T^n \times B \times J_0 \to \R$ in the $\omega$-Kolmogorov normal form.

On the other hand, we consider time-dependent vector fields on the torus $Z$ in such a way that
\begin{equation*}
Z:\T^n \times J_0 \to \R^n, \quad Z(q,t) = \omega + P(q,t) 
\end{equation*}
with $\omega \in \R^n$.

We introduce the following notation, which we will use in the rest of this work. 
For every functions $f$ defined on $\T^n \times B \times J_\upsilon$ and for fixed $t \in J_\upsilon$, let $f^t$ be the function defined on $\T^n \times B$ such that 
\begin{equation*}
f^t(q,p) = f(q,p,t)
\end{equation*}
for all $(q,p) \in \T^n \times B$. In addition, for fixed $p \in B$, let $f_{p}$ be the function defined on $\T^n \times J_\upsilon$ such that 
\begin{equation*}
f_{p}(q,t) = f(q,p,t)
\end{equation*}
for all $(q,t) \in \T^n \times J_\upsilon$. In accordance with the above notations, for fixed $(p,t) \in B \times J_\upsilon$, $f_p^t$ is the function defined on $\T^n$ such that
\begin{equation*}
f^t_{p}(q) = f(q,p,t)
\end{equation*}
for all $q \in \T^n$.

\subsection{Finitely Differentiable case}\label{IntroGD}

Here, we are interested in Hölder's classes of functions $C^\sigma$. We refer to Appendix \ref{A} for a very brief introduction. More specifically, in order to quantify the regularity of smooth functions, we introduce the following space. Given positive real parameters $\sigma \ge 0$ and $\upsilon \ge 0$, we have the following definition 
\begin{definition}
\label{SGD}
Let $\mathcal{S}^\upsilon_\sigma$ be the space of functions $f$ defined on  $\T^n \times B \times J_\upsilon $ such that $f \in C( \T^n \times B \times J_\upsilon)$  and, for all $t \in J_\upsilon$, $f^t \in C^\sigma( \T^n \times B)$.
\end{definition}
We use this notation also for functions defined on $\T^n \times J_\upsilon$, this will be specified by the context. Furthermore, for a positive integer $k \ge 0$, we have the following space of functions
\begin{definition}
\label{barSGD}
Let $\mathcal{\bar S}^\upsilon_{\sigma, k}$ be the space of functions $f$ such that 
\begin{equation*}
f \in \mathcal{S}^\upsilon_{\sigma + k} \hspace{3mm} \mbox{and} \hspace{3mm} \partial^i_{(q,p)}f \in \mathcal{S}^\upsilon_{\sigma + k-i}
\end{equation*}
for all $0 \le i \le k$.
\end{definition}
In the above definition, $\partial^i_{(q,p)}$ stands for the partial derivatives of order $i$ with respect to the variables $q$ and $p$. Conventionally $f = \partial^0_{(q,p)} f$. In other words, $f \in \mathcal{\bar S}^\upsilon_{\sigma, k}$ if $f \in \mathcal{S}^\upsilon_{\sigma + k}$ and  $\partial^i_{(q,p)}f \in C(\T^n \times B \times J_\upsilon)$ for all $0 \le i \le k$. That is, $f^t \in C^{\sigma + k}(\T^n \times B)$ for all $t \in J_\upsilon$ and the partial derivatives of $f$ with respect to $(q,p)$ are continuous until the order $k$. It is straightforward to verify that $\mathcal{\bar S}^\upsilon_{\sigma, 0} = \mathcal{S}^\upsilon_\sigma$. 

 In order to measure the decay in time of the perturbations, we introduce positive, decreasing, integrable functions $\mathbf{u}$ on $J_0$ and we denote 
\begin{equation*}
\mathbf{\bar u}(t) = \int_t^{+\infty} \mathbf{u}(\tau) d\tau
\end{equation*}
for all $t \in J_0$.

Now, we have everything we need to state the following theorem. Given $\omega \in \R^n$  and $\sigma \ge 1$, we consider a time-dependent Hamiltonian $H$ of the form 
\begin{equation}
\label{H1GD}
\begin{cases}
H : \T^n \times B \times J_0 \longrightarrow \R\\
H(q,p,t) = h(q,p,t) + f(q,p,t),\\
h \in \mathcal{K}_\omega\\
f_0, \left(\partial_p f\right)_0, \partial_p^2 H \in \mathcal{\bar S}^0_{\sigma, 2}\\
\displaystyle \sup_{t \in J_0}| f^t_0|_{C^{\sigma +2}} < \infty, \quad \sup_{t \in J_0}|\partial_p^2 H^t|_{C^{\sigma +2}} < \infty, \\
\displaystyle |\left(\partial_q f\right)^t_0|_{C^{\sigma +1}} \le \mathbf{a}(t), \quad \displaystyle |\left(\partial_p f\right)^t_0|_{C^{\sigma +2}}\le \mathbf{b}(t) \quad \mbox{for all $t \in J_0$},
\end{cases}
\tag{$*_A$}
\end{equation}
where $\mathbf{a}$, $\mathbf{b}$ are positive, decreasing, integrable functions on $J_0$.

We assume that there exists $\upsilon \ge 0$, such that $\mathbf{a}$ and $\mathbf{b}$ satisfy the following conditions
\begin{equation}
\label{propabGD}
\begin{cases}
\mathbf{\bar a}(t) \le \Lambda \mathbf{b}(t)\\
\mathbf{\bar a}(t) \mathbf{b}(t) \le \Lambda\mathbf{a}(t)\mathbf{\bar b}(t)
\end{cases}
\tag{$\#$}
\end{equation}
for all $t \in J_\upsilon$ and a suitable constant $\Lambda$. 

\begin{theoremx}
\label{Thm1GD}
Let $H$ be as in ~\eqref{H1GD} with $\mathbf{a}$ and $\mathbf{b}$ satisfying~\eqref{propabGD}. Then, there exist $\tilde h \in \mathcal{K}_\omega$ and a Lagrangian $C^\sigma$-asymptotic KAM torus $\varphi^t$ associated to $(X_H, X_{\tilde h}, \varphi_0)$.
\end{theoremx}

We begin with two examples of functions $\mathbf{a}$ and $\mathbf{b}$ satisfying~\eqref{propabGD}. First, we consider the case of exponential decay. Let $\mathbf{a}$ and $\mathbf{b}$ be the following functions
\begin{equation*}
\mathbf{a}(t) = e^{-\lambda_1 t}, \quad \mathbf{b}(t) = e^{-\lambda_2 t},
\end{equation*}
for some positive parameters $\lambda_1 \ge \lambda_2 >0$. It is straightforward to see that~\eqref{propabGD} is verified for all $t \in J_0$ and $\Lambda \ge \max\{{\lambda_2 \over \lambda_1}, {1 \over \lambda_1}\}$. 

The following example is more interesting than the previous one. It is about polynomial decay. We consider
\begin{equation*}
\mathbf{a}(t) = {1 \over t^{l+1}}, \quad \mathbf{b}(t) = {1 \over t^l},
\end{equation*}
for a positive real parameter $l >1$. This couple of functions satisfy~\eqref{propabGD} for all $t \in J_1$ with $\Lambda=1$.

The previous theorem proves the existence of a $C^\sigma$-asymptotic KAM torus $\varphi^t$ of the form
\begin{equation*}
\varphi^t(q) = (q + u^t(q), v^t(q))
\end{equation*}
for all $q \in \T^n$ and $t$ sufficiently large, where $\mathrm{id} + u^t$ is a diffeomorphism of the torus for all fixed $t$.  Furthermore, we also obtain some information about the decay in time of $u$ and $v$. More specifically
\begin{equation*}
|u^t|_{C^\sigma} \le C\mathbf{\bar b}(t), \quad |v^t|_{C^\sigma} \le C\mathbf{\bar a}(t),
\end{equation*}
for all $t$ large enough and for a suitable constant $C$.

Concerning time-dependent perturbations of constant vector fields on the torus, given $\sigma \ge 1$ and $\omega \in \R^n$, we consider the following time-dependent vector field
\begin{equation}
\label{ZGDIntro}
\begin{cases}
Z:\T^n \times J_0 \longrightarrow \R^n\\
Z(q,t) =  \omega + P(q,t)\\
P\in \mathcal{\bar S}^0_{\sigma,1},\\
|P^t|_{C^{\sigma+1}} \le \mathbf{P}(t) \quad \mbox{for all $t \in J_0$},
\end{cases}
\tag{$Z_A$}
\end{equation}
where $\mathbf{P}$ is a positive, decreasing, integrable function on $J_0$. 

\begin{corollaryx}
\label{Cor1GD}
Let $Z$ be as in~\eqref{ZGDIntro}. Then, there exists a $C^{\sigma}$-asymptotic KAM torus $\psi^t$ associated to $(Z, \omega, \mathrm{Id})$.
\end{corollaryx}

We observe that if $\mathbf{P}$ is not integrable on $J_0$, then, in general, it does not exist a $C^\sigma$-asymptotic KAM torus associated to $(Z, \omega, \mathrm{Id})$.  Let $\hat Z$ be a time-dependent vector field on $\T^1 \times J_0$ of the form
\begin{equation*}
\hat Z(q,t) = \omega + \hat P(t),
\end{equation*}
where $\omega \in \R$ and $\hat P(t)>0$ for all $t > 0$. We assume that
\begin{equation*}
\int_{t_0}^{+\infty} \hat P(\tau) d\tau = +\infty,
\end{equation*}
for all $t_0 \ge 0$. Let $\psi_{t_0, \hat Z}^t$ be the flow at time $t$ with initial time $t_0$ of $\hat Z$. Then, for all $q \in \T^n$ and  $t >0$ 
\begin{equation*}
\psi_{t_0, \hat Z}^{t_0 +t} (q) = q + \omega t + \int_{t_0}^{t_0 + t} \hat P(\tau) d\tau
\end{equation*}
and hence, for all $q \in \T^n$
\begin{eqnarray*}
\left| \psi_{t_0, \hat Z}^{t_0 +t} (q) - q - \omega t\right| &=&  \int_{t_0}^{t_0 + t} \hat P(\tau) d\tau.
\end{eqnarray*}
Therefore, taking the limit for $t \to +\infty$, the right-hand side of the latter diverges to $+\infty$. This means that does not exist an asymptotic quasiperiodic solution associated to $(\hat Z, \omega, \mathrm{Id})$ and thus does not exist a $C^\sigma$-asymptotic KAM torus associated to $(\hat Z, \omega, \mathrm{Id})$ (we refer to Proposition \ref{proplinkI}).

\subsection{Real analytic case}
As mentioned above, we state the real analytic version of the previous results. For some $s > 0$, we define complex domains 
\begin{equation*}
\T_s^n = \{q \in \C^n/\Z^n : | \operatorname{Im}(q)| \le s \}, \quad B_s =\{p \in \C^n : |p| \le s\},
\end{equation*}
and, given $\upsilon\ge 0$, we introduce the following space of functions.

\begin{definition}
\label{AsetfunctionGDA}
Let $\mathcal{A}^\upsilon_s$ be the space of the functions $f$ defined on  $\T^n_s \times B_s \times J_\upsilon$ such that $f \in C(\T^n_s \times B_s \times J_\upsilon)$ and, for all $t \in J_\upsilon$, $f^t$ is real analytic on $\T^n_s \times B_s$.
\end{definition}
We refer to the same notation for maps defined on $\T^n_s \times J_\upsilon$. For all  $k \in \Z^{2n}$ with $|k|\ge 1$, we let 
\begin{equation*}
\partial_{(q,p)}^{k} =  \partial^{k_1}_{q_1}...\partial^{k_n}_{q_n}\partial^{k_{n+1}}_{p_1}...\partial^{k_{2n}}_{p_n}
\end{equation*}
where $|k| = |k_1|+...+|k_{2n}|$. The following proposition concerns an important property of each $f \in \mathcal{A}^\upsilon_s$ that we will widely use in the rest of this work. 
\begin{proposition}
\label{propAAGD}
Let $f \in \mathcal{A}^\upsilon_s$, then, for all $k \in \Z^{2n}$ with $|k|\ge 1$,  $\partial_{(q,p)}^k f \in \mathcal{A}^\upsilon_{s'}$ for all $0 < s' <s$. 
\end{proposition}
\begin{proof}
For all $t \in J_\upsilon$, $\partial_{(q,p)}^{k}  f^t$ is real analytic on $\T^n_s \times B_s$ and hence on $\T^n_{s'} \times B_{s'}$. It remains to prove that $\partial_{(q,p)}^{k}  f \in C(\T^n_{s'} \times B_{s'} \times J_\upsilon)$. For all $(q_1,p_1, t_1)$, $(q_1,p_1, t_1) \in \T^n_{s'} \times B_{s'} \times J_\upsilon$, by Cauchy's inequality
\begin{eqnarray*}
|\partial_{(q,p)}^{k}  f(q_1,p_1, t_1) - \partial_{(q,p)}^{k}  f(q_2,p_2, t_2)| &\le& |\partial_{(q,p)}^{k}  f(q_1,p_1, t_1) - \partial_{(q,p)}^{k}  f(q_1,p_1, t_2)|\\
&+& |\partial_{(q,p)}^{k}  f(q_1,p_1, t_2) - \partial_{(q,p)}^{k}  f(q_2,p_2, t_2)|\\
&\le&{k_1!...k_{2n}! \over (s - s')^{|k|}} | f^{t_1} - f^{t_2}|_s\\
&+& \left|\partial_{(q,p)}^{k}  f^{t_2}(q_1,p_1) - \partial_{(q,p)}^{k}  f^{t_2}(q_2,p_2)\right|
\end{eqnarray*}
and hence, by the continuity of $f$ with respect to $t$ and the continuity of $\partial_{(q,p)}^{k}  f$ with respect to $(q,p)$, we have the claim.
\end{proof}
Given $\omega \in \R^n$ and a positive real parameter $s_0 >0$, we consider the following time-dependent Hamiltonian $H$ 
\begin{equation}
\label{H1AGD}
\begin{cases}
H : \T^n \times B \times J_0 \longrightarrow \R\\
H(q,p,t) = h(q,p,t) + f(q,p,t),\\
h \in \mathcal{K}_\omega\\
h ,f \in \mathcal{A}^0_{s_0}\\
\sup_{t \in J_0} |f^t_0|_{s_0} <\infty, \quad \sup_{t \in J_0}|\partial_p^2 H^t|_{s_0} <\infty\\
|\left(\partial_q f\right)^t_0|_{s_0}  \le \mathbf{a}(t), \quad |\left(\partial_p f\right)^t_0|_{s_0}  \le \mathbf{b}(t), \quad \mbox{for all $t \in J_0$}
\end{cases}
\tag{$*_B$}
\end{equation}
where $\mathbf{a}$, $\mathbf{b}$ are positive, decreasing, integrable functions on $J_0$.

\begin{theoremx}
\label{Thm1AGD}
Let $H$ be as in~\eqref{H1AGD} with $\mathbf{a}$ and $\mathbf{b}$ satisfying~\eqref{propabGD}. Then, there exist $\tilde h \in \mathcal{K}_\omega$ and a Lagrangian analytic asymptotic KAM torus $\varphi^t$ associated to $(X_H, X_{\tilde h}, \varphi_0)$.
\end{theoremx}

Similarly to Theorem \ref{Thm1GD}, we prove the existence of an analytic asymptotic KAM torus $\varphi^t$ of the form
\begin{equation*}
\varphi^t(q) = (q + u^t(q), v^t(q))
\end{equation*}
for all $q \in \T^n$ and $t$ sufficiently large, where $\mathrm{id} + u^t$ is a diffeomorphism of the torus for all fixed $t$. Furthermore, for a suitable constant $C$
\begin{equation*}
|u^t|_{s \over 4} \le C\mathbf{\bar b}(t), \quad |v^t|_{s \over 4} \le C\mathbf{\bar a}(t),
\end{equation*}
for all $t$ large enough. The proof of this theorem is essentially the same as that of Theorem \ref{Thm1GD} with suitable minor modifications.

Also in this case, we prove an analogous result regarding real analytic time-dependent perturbations of constant vector fields on the torus. Let $Z$ be a non-autonomous vector field on $\T^n \times J_0$ of the form
\begin{equation}
\label{ZaIntro}
\begin{cases}
Z:\T^n \times J_0 \longrightarrow \R^n,\\
Z(q,t) =  \omega + P(q,t)\\
P\in \mathcal{A}^0_{s_0},\\
|P^t|_{s_0} \le \mathbf{P}(t) \quad \mbox{for all $t \in J_0$}
\end{cases}
\tag{$Z_B$}
\end{equation}
where $ \omega \in \R^n$ and $0 < s_0 < 1$. We assume that $\mathbf{P}$ is a positive, decreasing, integrable function on $J_0$. 

\begin{corollaryx}
\label{Cor1AGD}
Let $Z$ be as in~\eqref{ZaIntro}. Then, there exists an analytic asymptotic KAM torus $\psi^t$ associated to $(Z, \omega, \mathrm{Id})$.
\end{corollaryx}
\begin{proof}
The proof is a straightforward application of Theorem \ref{Thm1AGD}. We consider the Hamiltonian $H$ defined on $\T^n \times B \times J_0$ of the form
\begin{equation*}
H(q, p, t) =  \omega \cdot p + P(q, t) \cdot p.
\end{equation*}
The latter satisfies the hypotheses of Theorem \ref{Thm1AGD}. Then, there exist $\tilde h \in \mathcal{K}_\omega$ and an analytic asymptotic KAM torus $\varphi^t$ associated to $(X_H, X_{\tilde h}, \varphi_0)$, where $\varphi_0$ is the trivial embedding previously introduced. Moreover, $\varphi^t = (\mathrm{id} + u^t, v^t)$ and, for all fixed $t$, $\mathrm{id} + u^t$ is a diffeomorphism of the torus. This concludes the proof of this theorem with $\psi^t = \mathrm{id} + u^t$.
\end{proof}

We point out that the previous argument does not work for Corollary \ref{Cor1GD} without asking $P\in \mathcal{\bar S}^0_{\sigma,2}$, hence a more regular perturbation. For this reason, the proof of Corollary  \ref{Cor1GD} is given in Section \ref{ProofC1GD}. 


\section{Asymptotic KAM tori}\label{AKAMtori}

We recall the definition of $C^\sigma$-asymptotic KAM torus. Let $\mathcal{P}$ be equal to $\T^n \times B$ or $\T^n$.
Given $\sigma \ge 0$ and a positive integer $k \ge 0$, we consider time-dependent vector fields $X^t$ and $X^t_0$ of class $C^{\sigma + k}$ on $\mathcal{P}$, for all $t \in J_{\upsilon}$, and an embedding $\varphi_0$ from $\T^n$ to $\mathcal{P}$ of class $C^\sigma$ such that
\begin{align}
\label{hyp1introKAM}
& \displaystyle \lim_{t \to +\infty} |X^t - X^t_0|_{C^{\sigma + k}} = 0,\\
\label{hyp2introKAM}
& X_0 (\varphi_0(q), t) =\partial_q \varphi_0(q)\omega \hspace{2mm} \mbox{for all $(q, t) \in \T^n \times J_\upsilon$,}
\end{align}
with $\omega \in \R^n$.
\begin{definitionnull}[Definition \ref{asymKAMtoriHol}]
We assume that $(X, X_0, \varphi_0)$ satisfy~\eqref{hyp1introKAM} and~\eqref{hyp2introKAM}. A family of $C^\sigma$ embeddings $\varphi^t: \T^n \to \mathcal{P}$ is a $C^\sigma$-asymptotic KAM torus associated to $(X, X_0, \varphi_0)$ if there exists $\upsilon' \ge \upsilon \ge 0$ such that  
\begin{align}
\label{Hyp1AKAMT}
&  \lim_{t \to +\infty}   |\varphi^t - \varphi_0|_{C^\sigma} = 0,\\
\label{Hyp2AKAMT}
&  X (\varphi(q,t), t) = \partial_q \varphi(q, t) \omega + \partial_t \varphi(q, t), 
\end{align}
for all $(q, t) \in \T^n \times J_{\upsilon'}$. When $\mathcal{P}$ is a symplectic manifold with $\mathrm{dim} \mathcal{P} = 2n$, then we say that $\varphi^t$ is Lagrangian if $\varphi^t(\T^n)$ is Lagrangian for all $t$.  
\end{definitionnull}

In what follows, we analyze a series of properties of $C^\sigma$-asymptotic KAM tori. Similar results are obviously true in the case of analytic asymptotic KAM tori. 

We observe that we can rewrite~\eqref{Hyp2AKAMT} in terms of the flow of $X$. In fact, let $\psi^t_{t_0,X}$ be the flow at time $t$ with initial time $t_0$ of $X$ and $\psi^t_{t_0, \omega}(q) = q + \omega(t -t_0)$ for all $q \in \T^n$ and $t$, $t_0 \in J_{\upsilon'}$.

\begin{proposition}
\label{IntroP1}
If the flow $\psi^t_{t_0,X}$ is defined for all $t$, $t_0 \in J_{\upsilon'}$, then~\eqref{Hyp2AKAMT} is equivalent to 
\begin{equation}
\label{hyp1aKtIbiss}
\psi^t_{t_0, X} \circ \varphi^{t_0}(q) = \varphi^t \circ \psi^t_{t_0, \omega}(q),
\end{equation} 
for all $t$, $t_0 \in J_{\upsilon'}$ and $q \in \T^n$.
\end{proposition}
\begin{proof}
In this proof, we denote the time dependence by indexes. We assume~\eqref{Hyp2AKAMT} and we prove~\eqref{hyp1aKtIbiss}. 
For fixed $t_0$, let $\varphi^{-t_0}$ be the inverse map of $\varphi^{t_0}$. It suffices to show that $\psi^t_{t_0, X}$ and $\varphi^t \circ \psi_{t_0, \omega}^t \circ \varphi^{-t_0}$ verify the same differential equation. For all $x \in \varphi^{t_0}(\T^n)$
\begin{eqnarray*}
{d \over dt} \left(\varphi^t \circ \psi_{t_0, \omega}^t \circ \varphi^{-t_0}(x)\right) &=& \partial_q \varphi^t \left( \psi_{t_0, \omega}^t \left( \varphi^{-t_0}(x) \right)\right)\dot \psi_{t_0, \omega}^t \left( \varphi^{-t_0}(x)\right)\\
&+& \partial_t \varphi^t \left( \psi_{t_0, \omega}^t \left( \varphi^{-t_0}(x) \right)\right)\\
&=& \partial_q \varphi^t \left( \psi_{t_0, \omega}^t \left( \varphi^{-t_0}(x) \right)\right)\omega + \partial_t \varphi^t \left( \psi_{t_0, \omega}^t \left( \varphi^{-t_0}(x) \right)\right)\\
&=& X^t\circ \varphi^t \circ \psi_{t_0, \omega}^t \circ \varphi^{-t_0}(x),
\end{eqnarray*}
where $\dot \psi^t_{t_0,\omega}$ stands for the derivative with respect to $t$ of $\psi^t_{t_0,\omega}$, it is obviously equal to $\omega$. The last equality is a consequence of~\eqref{Hyp2AKAMT}. This concludes the first part of the proof. 

Now, we assume~\eqref{hyp1aKtIbiss} and we prove~\eqref{Hyp2AKAMT}. We fix $t_0 \in J_{\upsilon'}$, for all $t \in J_{\upsilon'}$ and $x \in \varphi^{t_0}(\T^n)$
\begin{eqnarray*}
{d \over dt} \left(\varphi^t \circ \psi_{t_0, \omega}^t \circ \varphi^{-t_0}(x)\right) &=& \dot \psi^t_{t_0, X}(x) = X^t \circ \psi^t_{t_0, X}(x) = X^t \circ \varphi^t \circ \psi_{t_0, \omega}^t \circ \varphi^{-t_0}(x).
\end{eqnarray*}
On the other hand, by the chain rule 
\begin{eqnarray*}
{d \over dt} \left(\varphi^t \circ \psi_{t_0, \omega}^t \circ \varphi^{-t_0}(x)\right) = \partial_q \varphi^t \left( \psi_{t_0, \omega}^t \left( \varphi^{-t_0}(x) \right)\right)\omega + \partial_t \varphi^t \left( \psi_{t_0, \omega}^t \left( \varphi^{-t_0}(x) \right)\right). 
\end{eqnarray*}
We know that $\varphi^{t_0}$ is an embedding, then there exists $q\in \T^n$ such that $\varphi^{t_0}(q) = x$. Thanks to the above equations
\begin{equation*}
X^t \circ \varphi^t \circ \psi_{t_0, \omega}^t (q) = \partial_q \varphi^t \left( \psi_{t_0, \omega}^t(q)\right)\omega + \partial_t \varphi^t \left( \psi_{t_0, \omega}^t (q)\right)
\end{equation*}
for all  $q \in \T^n$ and for all $t \in J_{\upsilon'}$. Letting $t = t_0$ we have the claim.  
\end{proof}

By the latter, it is straightforward to see that~\eqref{Hyp2AKAMT} is trivial. 

\begin{proposition}
If $\psi^t_{t_0,X}$ is defined for all $t$, $t_0 \in J_{\upsilon'}$, it is always possible to find a family of embeddings satisfying~\eqref{Hyp2AKAMT}
\end{proposition}
\begin{proof}
Let $\hat \varphi : \T^n \to \mathcal{P}$ be an embedding then, for all $t$, $t_0 \in J_{\upsilon'}$ and $q \in \T^n$, we consider
\begin{equation*}
\varphi^t(q) = \psi_{t_0, X}^t \circ \hat \varphi (q - \omega (t - t_0)).
\end{equation*} 
The latter is a family of embeddings satisfying~\eqref{hyp1aKtIbiss}. Indeed, by the above definition of $\varphi^t$ we have that $\varphi^{t_0}(q) = \hat \varphi(q)$ for all $q \in \T^n$. Then, by construction, $\varphi^t$ satisfies~\eqref{hyp1aKtIbiss} and thus~\eqref{Hyp2AKAMT}.
\end{proof}

Another important consequence of~\eqref{hyp1aKtIbiss} is the following.  

\begin{proposition}
\label{IntroP3}
We assume that $\psi^t_{t_0,X}$ is defined for all $t$, $t_0 \in \R$. If there exists a $C^\sigma$-asymptotic KAM torus $\varphi^t$ defined for all $t \ge \upsilon'$, then we can extend the set of definition for all $t \in \R$.
\end{proposition}
\begin{proof}
For all $q \in \T^n$, we consider
\begin{equation}
\label{Introphiextended}
\phi^t(q) = \begin{cases} \varphi^t(q) \hspace{4mm} \mbox{for all $t \ge \upsilon'$}\\
\psi_{\upsilon', X}^t \circ \varphi^{\upsilon'}(q - \omega(t - \upsilon'))  \hspace{4mm} \mbox{for all $t \le \upsilon'$}.\end{cases}
\end{equation}
This is a family of embeddings which verify~\eqref{Hyp1AKAMT} and~\eqref{Hyp2AKAMT}.
\end{proof}

Unfortunately, we can not deduce any asymptotic information for the family of embeddings~\eqref{Introphiextended} when $t \to -\infty$. 

Now, concerning the dynamics associated to a $C^\sigma$-asymptotic KAM torus, we recall the definition of asymptotically quasiperiodic solution and discuss some properties of these motions. 
\begin{definitionnull}[Definition \ref{defqpsolI}]
We assume that $(X, X_0, \varphi_0)$ satisfy~\eqref{hyp1introAKAM} and~\eqref{hyp2introAKAM}. An integral curve $g(t)$ of $X$ is an asymptotically quasiperiodic solution associated to $(X, X_0, \varphi_0)$ if  there exists $q \in \T^n$ in such a way that 
\begin{equation*}
\lim_{t \to +\infty} |g(t) - \varphi_0 (q + \omega (t -t_0))|=0. 
\end{equation*}
\end{definitionnull}

The following proposition proves that if $\varphi^t$ is a $C^\sigma$-asymptotic KAM torus associated to $(X, X_0, \varphi_0)$, then each initial point $\varphi^{t_0}(q)$ gives rise to an asymptotically quasiperiodic solution associated to $(X, X_0, \varphi_0)$.
\begin{proposition}
\label{proplinkI}
Let $\varphi^t$ be a $C^\sigma$-asymptotic KAM torus associated to $(X, X_0, \varphi_0)$. 
Then, for all $q \in \T^n $ and $t_0 \in J_{\upsilon}$, 
\begin{equation*}
g(t) = \psi^t_{t_0, X}\circ \varphi^{t_0}(q)
\end{equation*}
is an asymptotically quasiperiodic solution associated to $(X, X_0, \varphi_0)$. 
\end{proposition}
\begin{proof}
Thanks to~\eqref{hyp1aKtIbiss}
\begin{equation*}
g(t) = \psi^t_{t_0, X}\circ \varphi^{t_0}(q) = \varphi^{t}(q + \omega(t-t_0))
\end{equation*}
and hence, by~\eqref{Hyp1AKAMT}, we have the claim. 
\end{proof}

We conclude this section with an important property concerning the case when $X$ and $X_0$ are Hamiltonian vector fields. Let $\mathcal{P} = \T^n \times B$ and we assume $\varphi^t$ to be a $C^\sigma$-asymptotic KAM torus associated to $(X, X_0, \varphi_0)$. In the particular case of Hamiltonian systems, if the autonomous KAM torus $\varphi_0$ is Lagrangian, then $\varphi^t$ is Lagrangian for all $t$. M.Canadell and R. de la Llave prove it in the discrete case. Here, we prove it in the continuous case.
\begin{proposition}
\label{IntroP2}
Let $\varphi^t$ be a $C^\sigma$-asymptotic KAM torus associated to $(X, X_0,\varphi_0)$. If $\varphi_0$ is Lagrangian, then $\varphi^t$ is Lagrangian for all $t \in J_{\upsilon'}$. 
\end{proposition}
\begin{proof}
Let $\alpha = dp \wedge dq$ be the standard symplectic form on $\T^n  \times B$. For all fixed $t$, $t_0 \in J_{\upsilon'}$, the map $\psi^t_{t_0, X}$  is a symplectomorphism. This means that $(\psi_{t_0, X}^t)^*\alpha = \alpha$  for all fixed $t$, $t_0 \in J_{\upsilon'}$. By~\eqref{hyp1aKtIbiss}, for all $t_0 \in J_{\upsilon'}$ and $t\ge 0$
\begin{equation}
\label{hyp1bissbissGD}
\psi_{t_0, X}^{t_0 + t} \circ \varphi^{t_0} = \varphi^{t_0 + t} \circ \psi_{t_0, \omega}^{t_0+t}
\end{equation}
and, taking the pull-back with respect to the standard form $\alpha$ on both sides of the latter, we obtain
\begin{equation*}
(\varphi^{t_0})^*(\psi_{t_0, X}^{t_0 + t})^* \alpha =  ( \psi_{t_0, \omega}^{t_0+t})^*(\varphi^{t_0 + t})^* \alpha.
\end{equation*}
We know that $\psi_{t_0, X}^{t_0 + t}$ is symplectic then, replacing $(\psi_{t_0, X}^{t_0 + t})^* \alpha = \alpha$ on the left hand side of the above equation, we have
\begin{equation*}
(\varphi^{t_0})^* \alpha =  ( \psi_{t_0, \omega}^{t_0+t})^*(\varphi^{t_0 + t})^* \alpha.
\end{equation*}
We want to prove that $\left((\varphi^{t_0})^* \alpha\right)_q = 0$ for all $q\in \T^n$, where $\left((\varphi^{t_0})^* \alpha\right)_q $ stands for the symplectic form calculated on $q\in \T^n$. The idea consists in verifying that, for all fixed $q\in \T^n$, the limit when $t \to +\infty$ on the right-hand side of the above equation converges to zero. Then, taking the limit for $t \to +\infty$ on both sides of the latter, we have the claim. 

We introduce the following notation
\begin{equation*}
\varphi^t(q) = (U^t(q), V^t(q)), \quad \varphi_0(q) = (U_0(q), V_0(q))
\end{equation*}
for suitable families of functions $U^t$, $V^t : \T^n \to \R^n$ and $U_0$, $V_0 : \T^n \to \R^n$. One can see that, for all $q \in \T^n$
\begin{equation*}
\left((\psi_{t_0, \omega}^{t_0+t})^*(\varphi^{t_0 + t})^* \alpha\right)_q = \sum_{1 \le i <j \le n} \alpha^{t_0 + t}_{ij}(q) dq_i \wedge dq_j,
\end{equation*}
where for all $1 \le i <j \le n$
\begin{eqnarray*}
\alpha^{t_0 + t}_{i,j}(q)  &=& \left(\partial_{q_i}V^{t_0 + t} \cdot \partial_{q_j}U^{t_0 + t}  - \partial_{q_j}V^{t_0 + t}  \cdot \partial_{q_i}U^{t_0 + t}  \right) \circ \psi_{t_0, \omega}^{t_0+t}(q).
\end{eqnarray*}
We observe that, for all $q \in \T^n$ and for all fixed $1 \le i <j \le n$
\begin{eqnarray*}
\partial_{q_i}V_0(q) \cdot \partial_{q_j}U_0(q) - \partial_{q_j}V_0(q) \cdot \partial_{q_i}U_0(q) =0
\end{eqnarray*}
because $\varphi_0$ is Lagrangian. Then,
\begin{eqnarray*}
\left|\alpha^{t_0 + t} _{i,j}\right|_{C^0}  &\le& \left|\partial_{q_i}V^{t_0 + t}  \cdot \partial_{q_j}U^{t_0 + t}  - \partial_{q_j}V^{t_0 + t}  \cdot \partial_{q_i}U^{t_0 + t}  \right|_{C^0}\\
&=& \left|\left(\partial_{q_i}V^{t_0 + t}  \cdot \partial_{q_j}U^{t_0 + t}  - \partial_{q_j}V^{t_0 + t}  \cdot \partial_{q_i}U^{t_0 + t} \right) - \left(\partial_{q_i}V_0 \cdot \partial_{q_j}U_0 - \partial_{q_j}V_0 \cdot \partial_{q_i}U_0\right)  \right|_{C^0}\\
&\le& \left|\partial_{q_i}V^{t_0 + t}  \cdot \partial_{q_j}U^{t_0 + t}  - \partial_{q_i}V_0 \cdot \partial_{q_j}U_0\right|_{C^0} + \left|\partial_{q_j}V^{t_0 + t}  \cdot \partial_{q_i}U^{t_0 + t}  - \partial_{q_j}V_0 \cdot \partial_{q_i}U_0\right|_{C^0},
\end{eqnarray*}
and we can estimate each term in the last line of the latter by
\begin{eqnarray*}
\left|V^{t_0 + t} \right|_{C^1} \left|U^{t_0 + t}  -U_0 \right|_{C^1} + \left|U_0\right|_{C^1} \left|V^{t_0 + t}  -V_0 \right|_{C^1} 
\end{eqnarray*}
multiplied by a suitable constant $C$. This concludes the proof of this proposition because, by~\eqref{Hyp1AKAMT}, the latter converges to $0$ if $t \to +\infty$. 
\end{proof}

\section{Proof of Theorem \ref{Thm1GD}}\label{ProofT1GD}
This section is devoted to the proof of Theorem \ref{Thm1GD}. To this end, we expand the Hamiltonian $H$ in~\eqref{H1GD} in a small neighbourhood of $0 \in B$, 
\begin{eqnarray*}
h(q,p,t) &=& h(q,0,t) + \partial_p h(q,0,t) \cdot p +  \int_0^1 (1 -\tau) \partial^2_p h(q, \tau p,t)d\tau \cdot p^2\\
f(q,p,t) &=& f(q, 0, t) + \partial_p f(q, 0,t) \cdot p +  \int_0^1 (1 -\tau)\partial^2_p f(q, \tau p, t) d\tau\cdot p^2,
\end{eqnarray*}
we can assume without loss of generality that $h(q,0,t)=0$ for all $(q,t) \in \T^n \times J_0$. Letting
\begin{eqnarray*}
\omega &=& \partial_p h(q,0,t)\\
a(q,t) &=& f(q,0,t)\\
b(q,t) &=& \partial_p f(q,0,t)\\
m(q,p,t) &=& \int_0^1 (1 -\tau) \left(\partial_p^2 h(q,\tau p, t) + \partial_p^2 f(q, \tau p, t)\right) d\tau \\
&=& \int_0^1 (1 -\tau) \partial_p^2 H(q,\tau p, t) d\tau,
\end{eqnarray*}
for a positive real parameter $\Upsilon \ge 1$, we can rewrite the Hamiltonian $H$ in the following form
\begin{equation}
\label{H2GD}
\begin{cases}
H : \T^n \times B \times J_0 \longrightarrow \R\\
H(q,p,t) = \omega \cdot p + a(q,t) + b(q,t) \cdot p + m(q,p,t) \cdot p^2,\\
a, b, \partial^2_p H \in \mathcal{\bar S}^0_{\sigma,2},\\
\sup_{t \in J_0} |a^t|_{C^{\sigma+2}} < \infty, \quad \sup_{t \in J_0}|\partial_p^2 H^t|_{C^{\sigma +2}} \le \Upsilon,\\
|\partial_q a^t|_{C^{\sigma+1}} \le \mathbf{a}(t), \quad |b^t|_{C^{\sigma+2}}  \le \mathbf{b}(t), \quad \mbox{for all $t \in J_0$}
\end{cases}
\tag{$**_A$}
\end{equation}
where $\mathbf{a}(t)$ and $\mathbf{b}(t)$ are the functions introduced in~\eqref{H1GD} satisfying~\eqref{propabGD}. This Hamiltonian is our new starting point. Furthermore, let $\tilde h$ be the following Hamiltonian 
\begin{equation*}
\tilde h(q,p,t) = h(q,p,t) + \int_0^1 (1 -\tau)\partial^2_p f(q, \tau p, t) d\tau\cdot p^2
\end{equation*}
for all $(q,p,t) \in \T^n \times B \times J_0$. Obviously $\tilde h \in \mathcal{K}_\omega$. Moreover, $X_H$ and $X_{\tilde h}$ verify~\eqref{hyp1introKAM}.

\subsection{Outline of the Proof of Theorem \ref{Thm1GD}}\label{OPGD}

We are looking for a $C^\sigma$-asymptotic KAM torus $\varphi^t$ associated to $(X_H, X_{\tilde h}, \varphi_0)$, where $H$ is the Hamiltonian in~\eqref{H2GD}, $\tilde h$ is the Hamiltonian previously defined and $\varphi_0$ is the trivial embedding $\varphi_0 : \T^n  \to \T^n \times B$, $\varphi_0(q) = (q,0)$. More specifically, for given $H$, we are searching for $\upsilon' \ge 0$ sufficiently large and suitable functions $u$, $v : \T^n \times J_{\upsilon'} \to \R^n$ such that 
\begin{equation*}
\varphi(q,t) = (q + u(q,t), v(q,t))
\end{equation*}
and in such a way that  $\varphi$, $u$ and $v$ satisfy the following conditions
\begin{align}
\label{hyp1GD2}
&   X_H(\varphi(q, t), t) -  \partial_q \varphi(q, t) \omega - \partial_t \varphi(q, t) = 0,\\
\label{hyp2GD2}
&  \lim_{t \to +\infty}  |u^t|_{C^\sigma} = 0, \quad \lim_{t \to +\infty}  |v^t|_{C^\sigma} = 0,
\end{align}
for all $(q,t) \in \T^n \times J_{\upsilon'}$. The parameter $\upsilon'$ is free and it will be chosen large enough in Lemma \ref{lemmautilethmGD} below.

The proof rests on the implicit function theorem. To this end, we need to introduce a suitable functional $\mathcal{F}$ given by~\eqref{hyp1GD2}. We consider
\begin{equation*}
\bar m(q,p,t) p = \left(\int_0^1  \partial_p^2 H(q,\tau p, t) d\tau \right) p = \partial_p \Big(m(q,p,t) \cdot p^2 \Big).
\end{equation*}
This is well defined because
\begin{eqnarray*}
\partial_p \Big(m(q,p,t) \cdot p^2 \Big) &=& \partial_p \left(\int_0^1 (1-\tau) \partial_p^2 H(q,\tau p, t) d\tau \cdot p^2\right)\\
&=& \partial_p \left(\int_0^1 (p - \xi) \partial_p^2 H(q,\xi, t) d\xi\right)\\
&=& \int_0^1  \partial_p^2 H(q,\xi, t) d\xi = \left(\int_0^1  \partial_p^2 H(q,\tau p, t) d\tau \right) p,
\end{eqnarray*}
where the second equality of the latter is due to the change of variables $\xi = \tau p$. Going back to the definition of the functional $\mathcal{F}$, we observe that the Hamiltonian system associated to the Hamiltonian $H$  is equal to 
\begin{equation*}
X_H(q,p,t)  = \begin{pmatrix} \omega + b(q,t) + \bar m(q,p,t) p  \\
-\partial_q a(q, t) - \partial_q b(q,t)p - \partial_q m(q,p,t) p^2\end{pmatrix},
\end{equation*}
where we recall that $H$ is the Hamiltonian defined by~\eqref{H2GD}. We introduce
\begin{equation*}
\tilde \varphi(q,t) = (q + u(q,t), v(q,t), t), \quad \tilde u(q,t) = (q + u(q,t), t),
\end{equation*}
for all $(q,t) \in \T^n \times J_{\upsilon'}$. Composing the Hamiltonian system $X_H$ with $\tilde \varphi$, we can write $X_H \circ \tilde \varphi $ in the following form
\begin{small}
\begin{equation*}
X_H \circ \tilde \varphi (q, t) = \begin{pmatrix}\omega + b \circ \tilde u (q, t) + \bar m \circ \tilde \varphi  (q, t) v (q, t)  \\
-\partial_q a\circ \tilde u (q, t) - \partial_q b \circ \tilde u (q, t) v (q, t) - \partial_q m \circ \tilde \varphi (q, t) \cdot v (q, t)^2\end{pmatrix}
\end{equation*}
\end{small}
for all $(q,t) \in \T^n \times J_{\upsilon'}$ and moreover
\begin{equation*}
\partial_q \varphi(q, t)\omega + \partial_t \varphi(q, t) = \begin{pmatrix} \omega + \partial_q u(q, t)\omega + \partial_t u (q, t) \\
\partial_q v(q, t)\omega + \partial_t v (q, t)\end{pmatrix}
\end{equation*}
for all $(q,t) \in \T^n \times J_{\upsilon'}$. We define
\begin{equation*}
\nabla u (q,t)\Omega = \partial_q u(q,t) \omega + \partial_t u(q,t), \quad \nabla v(q,t) \Omega = \partial_q v(q,t) \omega + \partial_t v(q,t)
\end{equation*}
for all $(q, t) \in \T^n \times J_{\upsilon'}$. Then, we can rewrite~\eqref{hyp1GD2} in the following form 
\begin{eqnarray}
\label{F.0GD}
\begin{pmatrix}b \circ \tilde u + \left(\bar m \circ \tilde \varphi  \right) v - \left(\nabla u \right)\Omega \\
-\partial_q a\circ \tilde u  - \left(\partial_q b \circ \tilde u \right) v  - \left(\partial_q m \circ \tilde  \varphi   \right) \cdot v^2 - \left(\nabla v \right)\Omega 
\end{pmatrix} = \begin{pmatrix} 0\\0\end{pmatrix}.
\end{eqnarray}
This is composed of sums and products of functions defined on $(q, t) \in \T^n \times J_{\upsilon'}$, we have omitted the arguments $(q, t)$ in order to achieve a more elegant form. We keep this notation for the rest of this proof. Over suitable Banach spaces, that we will specify later, let $\mathcal{F}$ be the following functional
\begin{equation*}
\mathcal{F}(a, b, m, \bar m,u, v) = (F_1(b, \bar m,u,v), F_2(a,b, m,u,v))
\end{equation*}
with 
\begin{eqnarray*}
F_1(b, \bar m,u,v) &=& b \circ \tilde u + \left(\bar m \circ \tilde \varphi \right) v - \left(\nabla u \right)\Omega,\\
F_2(a,b, m,u,v) &=& \partial_q a\circ \tilde u  + \left(\partial_q b \circ \tilde u \right) v  + \left(\partial_q m \circ \tilde \varphi \right) \cdot v^2 + \left(\nabla v \right)\Omega.
\end{eqnarray*}
The latter is obtained by~\eqref{F.0GD} and we observe that for all $m$ and $\bar m$, 
\begin{equation*}
\mathcal{F}(0, 0, m, \bar m,0, 0)=0.
\end{equation*}
We can reformulate our problem in the following form. For fixed $m$ and $\bar m$ in a suitable Banach space and for $(a,b)$ sufficiently close to $(0,0)$, we are looking for some functions $u$, $v$ in such a way that $\mathcal{F}(a, b, m, \bar m,u, v) = 0$ and the asymptotic conditions~\eqref{hyp2GD2} are satisfied. 

Concerning the associated linearized problem, the differential of $\mathcal{F}$ with respect to the variables $(u,v)$ calculated on $(0,0,m, \bar m, 0,0)$ is equal to
\begin{equation*}
D_{(u,v)} \mathcal{F}(0,0,m,\bar m, 0, 0)(\hat u, \hat v) = (\bar m_0\hat v - \left(\nabla \hat u\right)\Omega,  \left(\nabla \hat v\right)\Omega)
\end{equation*}
where, in according to the notation previously introduced, for all $(q, t) \in \T^n \times J_{\upsilon'}$ we let $\bar m_0(q,t) = \bar m(q,0,t)$.

The proof of this theorem is a straightforward application of the implicit function theorem if we assume the following norm
\begin{equation*}
\sup_{t \in J_0} |a^t|_{C^{\sigma+2}} + \sup_{t \in J_0}{|\partial_q a^t|_{C^{\sigma+1}} \over \mathbf{a}(t)}, \hspace{7mm}   \sup_{t \in J_0}{|b^t|_{C^{\sigma+2}} \over \mathbf{b}(t)},
\end{equation*} 
to be sufficiently small. To avoid this smallness assumption, we study the problem from another point of view.
We are looking for a $C^\sigma$-asymptotic KAM torus defined for $t$ sufficiently large in such a way 
\begin{equation*}
|\partial_q a^t|_{C^{\sigma+1}} , \hspace{7mm}   |\partial_q b^t|_{C^{\sigma+1}} 
\end{equation*} 
are sufficiently small. It suffices for proving the existence of functions $(u,v)$ satisfying~\eqref{hyp1GD2} and~\eqref{hyp2GD2}.

The following four sections are devoted to the proof of Theorem \ref{Thm1GD}. In the first, we introduce suitable Banach spaces on which the previous functional is defined. The second is dedicated to solving the homological equation, which is the main tool to prove that $D_{(u,v)} \mathcal{F}(0,0,m,\bar m, 0, 0)$ is invertible. In the penultimate section, we verify that $\mathcal{F}$ is well-defined and satisfies the hypotheses of the implicit function theorem. Finally, the last section concludes the proof of this theorem.

\subsection{Preliminary Settings}\label{PrelSettGD}

Given positive real parameters $\sigma \ge 0$, $\upsilon \ge 0$ and a positive integer $k \ge 0$, we recall that $\mathcal{S}^\upsilon_\sigma$ and $\mathcal{\bar S}^\upsilon_{\sigma, k}$ are respectively the spaces of functions defined by Definition \ref{SGD} and Definition \ref{barSGD}.
We introduce the following norm that we will widely use in the rest of this section. For every $f \in \mathcal{S}^\upsilon_{\sigma}$ and for a positive real function $\mathbf{u}(t)$ defined on $J_\upsilon$, we define
\begin{equation}
\label{defnormGD}
|f|_{\sigma, \mathbf{u}}^\upsilon = \sup_{t \in J_\upsilon}{|f^t|_{C^{\sigma}} \over \mathbf{u}(t)}.
\end{equation}
Furthermore, we recall that for all positive, decreasing, integrable functions $\mathbf{u}$ on $J_\upsilon$, we let $\mathbf{\bar u}$ be  
\begin{equation*}
\mathbf{\bar u}(t) = \int_t^{+\infty} \mathbf{u}(\tau) d\tau
\end{equation*}
for all $t \in J_\upsilon$.

Now, let  $\sigma \ge 1$, $\upsilon \ge 0$ and $\Upsilon \ge 1$ be the positive parameters introduced in~\eqref{H2GD} and~\eqref{propabGD}. For $\upsilon' \ge \upsilon \ge0$ that will be chosen later, we consider the following Banach spaces $\left(\mathcal{A}, |\cdot |\right)$, $\left(\mathcal{B}, |\cdot |\right)$, $\left(\mathcal{U}, |\cdot |\right)$, $\left(\mathcal{V}, |\cdot |\right)$, $\left(\mathcal{Z}, |\cdot |\right)$ and $\left(\mathcal{G}, |\cdot |\right)$ (see Appendix \ref{BanSp})
\vspace{5mm}
\begin{eqnarray*}
\mathcal{A} &=& \Big\{a : \T^n \times J_{\upsilon'} \to \R  \hspace{1mm}| \hspace{1mm} a \in \mathcal{\bar S}^{\upsilon'}_{\sigma, 2} \hspace{1mm} \mbox{and} \hspace{1mm} |a| =|a|^{\upsilon'}_{\sigma+2, 1} + |\partial_q a|^{\upsilon'}_{\sigma+1, \mathbf{a}} < \infty\Big\}\\
\mathcal{B} &=& \Big\{b : \T^n \times J_{\upsilon'}  \to \R^n  \hspace{1mm}| \hspace{1mm} b \in \mathcal{\bar S}^{\upsilon'}_{\sigma, 2} , \hspace{1mm} \mbox{and} \hspace{1mm} |b| =|b|^{\upsilon'}_{\sigma+2, \mathbf{b}} < \infty\Big\}\\
\mathcal{U} &=& \Big\{u :\T^n \times J_{\upsilon'}  \to \R^n  \hspace{1mm}| \hspace{1mm} u,\left( \nabla u\right) \Omega \in \mathcal{S}^{\upsilon'}_{\sigma} \\
&& \mbox{and} \hspace{1mm} |u| = \max\{|u|^{\upsilon'}_{\sigma, \mathbf{\bar b}}, | \left(\nabla u \right) \Omega|^{\upsilon'}_{\sigma, \mathbf{b}} \} < \infty\Big\}\\
\mathcal{V} &=& \Big\{v :\T^n \times J_{\upsilon'}  \to \R^n  \hspace{1mm}| \hspace{1mm} v,\left( \nabla v\right) \Omega \in \mathcal{S}^{\upsilon'}_{\sigma}  \\
&& \mbox{and} \hspace{1mm} |v| = \max\{|v|^{\upsilon'}_{\sigma, \mathbf{\bar a}}, |\left( \nabla v \right) \Omega|^{\upsilon'}_{\sigma, \mathbf{a}} \} < \infty\Big\}\\
\mathcal{Z} &=& \Big\{z : \T^n \times J_{\upsilon'}  \to \R^n  \hspace{1mm}| \hspace{1mm} z \in \mathcal{S}^{\upsilon'}_{\sigma}   , \hspace{1mm} \mbox{and} \hspace{1mm} |z| =|z|^{\upsilon'}_{\sigma, \mathbf{b}} < \infty\Big\}\\\mathcal{G} &=& \Big\{g : \T^n \times J_{\upsilon'} \to \R  \hspace{1mm}| \hspace{1mm} g \in \mathcal{S}^{\upsilon'}_{\sigma}   \hspace{1mm} \mbox{and} \hspace{1mm} |g| =|g|^{\upsilon'}_{\sigma, \mathbf{a}} < \infty\Big\}\\
\end{eqnarray*} 
where, in the definition of $\mathcal{A}$, the norm $|a|^{\upsilon'}_{\sigma+2, 1} = \sup_{t \in J_{\upsilon'}}|a^t|_{C^{\sigma+2}}$. This means that $1$ stands for the function identically equal to $1$ for all $t \in J_{\upsilon'}$. Let $M_n$ be the set of the $n$-dimensional matrices. We introduce another Banach space $\left(\mathcal{M}, |\cdot |\right)$ in such a way that 
\begin{equation*}
\mathcal{M} = \Big\{m : \T^n \times B \times J_{\upsilon'} \to M_n  \hspace{1mm}| \hspace{1mm} m \in \mathcal{\bar S}^{\upsilon'}_{\sigma, 2} \hspace{1mm} \mbox{and} \hspace{1mm} |m| =|m|^{\upsilon'}_{\sigma+2, 1} \le \Upsilon\Big\}
\end{equation*}
where $\Upsilon$ is the positive parameter in~\eqref{H2GD}. Now, we have everything we need to define more precisely the functional $\mathcal{F}$ introduced in the previous section. Let $\mathcal{F}$ be the following functional
\begin{equation*}
\mathcal{F} :  \mathcal{A} \times \mathcal{B} \times \mathcal{M} \times \mathcal{M} \times  \mathcal{U} \times  \mathcal{V} \longrightarrow \mathcal{Z} \times \mathcal{G}\\
\end{equation*}
\begin{equation*}
\mathcal{F}(a, b,m , \bar m,  u, v) = (F_1(b, \bar m, u,v), F_2(a,b,m,u,v))
\end{equation*}
with 
\begin{eqnarray*}
F_1(b,\bar m ,u,v) &=& b \circ \tilde u + \left(\bar m \circ \tilde \varphi \right) v - \left(\nabla u \right)\Omega,\\
F_2(a,b, m,u,v) &=& \partial_q a\circ \tilde u  + \left(\partial_q b \circ \tilde u \right) v  + \left(\partial_q m \circ \tilde \varphi \right) \cdot v^2 + \left(\nabla v \right)\Omega.
\end{eqnarray*}

\subsection{Homological Equation}\label{SecHEGD}
Given $\sigma \ge 0$, $\upsilon \ge 0$ and $\omega \in \R^n$, in this section, we solve the following equation for the unknown $\varkappa : \T^n \times J_\upsilon \to \R$
\begin{equation}
\label{HEGD}
\begin{cases}
 \omega \cdot \partial_q \varkappa(q,t) + \partial_t \varkappa(q,t) = g(q,t),\\
g \in \mathcal{S}^{\upsilon}_{\sigma}  \quad |g|_{\sigma, \mathbf{g}}^\upsilon < \infty,
\end{cases}
\tag{$HE_A$}
\end{equation}
where $\mathbf{g}(t)$ is a positive, decreasing, integrable function on $J_\upsilon$ and $g : \T^n \times J_\upsilon \to \R$ is given. 

\begin{lemma}[\textbf{Homological Equation}]
\label{homoeqlemmaGD} 
There exists a unique solution $\varkappa \in  \mathcal{S}^{\upsilon}_{\sigma} $ of~\eqref{HEGD} such that 
\begin{equation}
\label{varkappaasym}
\lim_{t \to +\infty} |\varkappa^t|_{C^0} = 0.
\end{equation}
Moreover, 
\begin{equation*}
|\varkappa|^\upsilon_{\sigma, \mathbf{\bar g}} \le |g|_{\sigma, \mathbf{g}}^\upsilon.
\end{equation*}
\end{lemma}
\begin{proof}
\textit{Existence}: Let us define the following transformation 
\begin{equation*}
h: \T^n \times J_\upsilon \to \T^n \times J_\upsilon, \quad h(q,t) = (q -  \omega t , t),
\end{equation*}
which is the key to solving the homological equation.

We claim that it is enough to prove the first part of this lemma for the much simpler equation
\begin{equation}
\label{homoeqnewGD}
\partial_t \kappa = g(q+ \omega t, t).
\end{equation}
As a matter of fact, if $\kappa$ is a solution of the latter satisfying the asymptotic condition~\eqref{varkappaasym}, then $\chi = \kappa \circ h$ is a solution of (\ref{HEGD}) satisfying the same asymptotic condition and viceversa. For the sake of clarity, we prove this claim. Let $\varkappa$ be a solution of~\eqref{HEGD} verifying the asymptotic condition~\eqref{varkappaasym}, then 
\begin{equation*}
\partial_t (\varkappa \circ h^{-1}) = \partial_q \varkappa \circ h^{-1} \cdot \omega + \partial_t \varkappa \circ h^{-1} = g\circ h^{-1},
\end{equation*}
where the last equality is due to~\eqref{HEGD}. This implies that $\kappa = \varkappa \circ h^{-1}$ is a solution of~\eqref{homoeqnewGD} and by
\begin{equation*}
|\kappa^t|_{C^0} = |\left(\varkappa \circ h^{-1}\right)^t|_{C^0} \le |\varkappa^t|_{C^0}
\end{equation*}
$\kappa = \varkappa \circ h^{-1}$ satisfies the asymptotic condition because $\varkappa$ does. Viceversa, let $\kappa$ be a solution of~\eqref{homoeqnewGD} satisfying the asymptotic condition~\eqref{varkappaasym}, then
\begin{equation*}
\partial_q (\kappa \circ h)\cdot \omega  + \partial_t(\kappa \circ h) = \partial_q \kappa \circ h \cdot \omega - \partial_q \kappa \circ h\cdot \omega  + \partial_t \kappa \circ h = g.
\end{equation*}
By~\eqref{homoeqnewGD}, we have the last equality of the latter. Hence, $\kappa \circ h$ is a solution of~\eqref{HEGD}. Moreover, thanks to 
\begin{equation*}
|\varkappa^t|_{C^0} = |\left(\kappa \circ h\right)^t|_{C^0} \le |\kappa^t|_{C^0}
\end{equation*}
$\varkappa = \kappa \circ h$ satisfies the asymptotic condition~\eqref{varkappaasym}. This proves the claim.

For all $q \in \T^n$ a solution of~\eqref{homoeqnewGD} exists and 
\begin{equation*}
\kappa(q,t) = e(q) + \int_\upsilon^t g(q+ \omega \tau, \tau)d\tau
\end{equation*}
with a function $e$ defined on the torus. We have to choose $e$ in such a way that $\kappa$ satisfies the following asymptotic condition for all fixed $q \in \T^n$
\begin{equation*}
0 = \lim_{t \to +\infty} \kappa(q,t) = e(q) + \int_\upsilon^{+\infty} g(q+ \omega \tau, \tau)d\tau. 
\end{equation*}
There is only one possible choice for $e$, and that is
\begin{equation*}
e(q) = -\int_\upsilon^{+\infty} g(q+ \omega \tau, \tau)d\tau.
\end{equation*}
This implies that 
\begin{equation*}
\kappa(q,t) = -\int_t^{+\infty} g(q+ \omega\tau, \tau)d\tau
\end{equation*}
is the solution of (\ref{homoeqnewGD}) we are looking for. Therefore, $e$ is well defined, indeed
\begin{equation*}
\left|\int_\upsilon^{+\infty} g(q+ \omega \tau, \tau)d\tau \right| \le |g|_{\sigma, \mathbf{g}}^\upsilon \int_\upsilon^{+\infty} \mathbf{g}(\tau) d\tau = |g|_{\sigma, \mathbf{g}}^\upsilon \mathbf{\bar g}(\upsilon) < \infty.
\end{equation*}
Moreover, 
\begin{equation*}
|\kappa^t|_{C^0} \le  \int_t^{+\infty} |g^\tau|_{C^0}d\tau  \le |g|_{\sigma, \mathbf{g}}^\upsilon \int_t^{+\infty} \mathbf{g}(\tau) d\tau = |g|_{\sigma, \mathbf{g}}^\upsilon \mathbf{\bar g}(t),
\end{equation*}
since $\mathbf{\bar g}(t)$ converges to $0$ when $t \to +\infty$, taking the limit for $t \to +\infty$ on both sides of the latter, we have that $|\kappa^t|_{C^0} \to 0$ when $t \to +\infty$. 
This concludes the first part of the proof because 
\begin{equation*}
\varkappa(q,t) = \kappa \circ h(q,t) =-\int_t^{+\infty} g(q+ \omega(\tau - t), \tau)d\tau
\end{equation*}
is the unique solution of (\ref{HEGD}) verifying~\eqref{varkappaasym} that we are looking for.

\textit{Regularity and Estimates}: We observe that $g \in \mathcal{S}^{\upsilon}_{\sigma} $ implies $\kappa \in \mathcal{S}^{\upsilon}_{\sigma} $ and hence $\varkappa = \kappa \circ h \in\mathcal{S}^{\upsilon}_\sigma$. Moreover, for all fixed $t \in J_\upsilon$
\begin{equation*}
|\varkappa^t|_{C^\sigma} \le |g|_{\sigma, \mathbf{g}}^\upsilon \mathbf{\bar g}(t).
\end{equation*}
Multiplying both sides of the latter by ${1 \over \mathbf{\bar g}(t)}$ and taking the sup for all $t \in J_\upsilon$, we prove the second part of this lemma. 
\end{proof}

\subsection{Regularity of $\mathcal{F}$}\label{RegFGD}
We recall the  definition of the functional $\mathcal{F}$,
\begin{equation*}
\mathcal{F} :  \mathcal{A} \times \mathcal{B} \times \mathcal{M} \times \mathcal{M} \times  \mathcal{U} \times  \mathcal{V} \longrightarrow \mathcal{Z} \times \mathcal{G}\\
\end{equation*}
\begin{equation*}
\mathcal{F}(a, b,m , \bar m,  u, v) = (F_1(b, \bar m, u,v), F_2(a,b,m,u,v))
\end{equation*}
with 
\begin{eqnarray*}
F_1(b,\bar m ,u,v) &=& b \circ \tilde u + \left(\bar m \circ \tilde \varphi \right) v - \left(\nabla u \right)\Omega,\\
F_2(a,b, m,u,v) &=& \partial_q a\circ \tilde u  + \left(\partial_q b \circ \tilde u \right) v  + \left(\partial_q m \circ \tilde \varphi \right) \cdot v^2 + \left(\nabla v \right)\Omega.
\end{eqnarray*}
Here, we verify that $\mathcal{F}$ satisfies the hypotheses of the implicit function theorem. Using the properties in Proposition \ref{Holder} (see Appendix \ref{A}) and~\eqref{propabGD}, one can prove that $\mathcal{F}$ is well defined, continuous and differentiable with respect to the variables $(u,v)$. Let $D_{(u,v)} \mathcal{F}$ be the differential of $\mathcal{F}$ with respect to $(u,v)$, we have that 
\begin{eqnarray}
\label{DF1GD}
D_{(u,v)} F_1(b,\bar m,u,v)(\hat u, \hat  v) &=& D_u F_1(b,\bar m,u,v)\hat u  + D_v F_1(b,\bar m,u,v)\hat v  \nonumber \\
&=& \left(\partial_q b \circ \tilde u \right) \hat  u + v^T \left(\partial_q \bar m \circ \tilde  \varphi \right)\hat u + v^T \left( \partial_p \bar m \circ \tilde  \varphi \right)\hat v \nonumber\\
&+& \left(\bar m \circ \tilde  \varphi \right)\hat v - \left(\nabla \hat u \right) \Omega\\
\label{DF2GD}
D_{(u,v)} F_2(a, b,m,u,v) (\hat u, \hat v) &=& D_u F_2(a,b,m,u,v)\hat u  + D_v F_2(a,b,m,u,v)\hat v  \nonumber\\
&=& \left(\partial^2_q a \circ \tilde u \right) \hat u  + v^T \left( \partial^2_q b \circ \tilde u \right) \hat u+  (v^T)^2 \left( \partial^2_q m \circ \tilde  \varphi \right)\hat u\nonumber \\
&+& \left(\partial_q b \circ \tilde u \right) \hat v  +  (v^T)^2 \left(\partial^2_{pq} m \circ \tilde  \varphi \right)\hat v + 2 v^T \left(\partial_q m \circ \tilde  \varphi \right)\hat v\nonumber\\
&+&  \left(\nabla \hat v \right)\Omega,
\end{eqnarray}
where $T$ stands for the transpose of a vector and $D_u$, $D_v$ are respectively the differentials with respect to $u$ and $v$. Furthermore, one can show that $D_{(u,v)}\mathcal{F}$ is continuous. Now, we observe that this differential calculated on $(0,0,m,\bar m, 0,0)$ is equal to
\begin{equation}
\label{DFGD}
D_{(u,v)} \mathcal{F}(0,0,m,\bar m, 0, 0)(\hat u, \hat v) = (\bar m_0\hat v - \left(\nabla \hat u\right)\Omega,  \left(\nabla \hat v\right)\Omega).
\end{equation}
In the following lemma, we verify that the latter is invertible for all fixed $m$, $\bar m \in \mathcal{M}$. First, to avoid a flow of constants, let $C(\cdot)$ be constants depending on $n$ and the other parameters in brackets.  On the other hand, $C$ stands for constants depending only on $n$.

\begin{lemma}
\label{lemmainvGD}
For all $(z,g) \in \mathcal{Z} \times \mathcal{G}$ there exists a unique $(\hat u, \hat v) \in \mathcal{U} \times \mathcal{V}$ such that 
\begin{equation*}
D_{(u,v)}\mathcal{F}(0,0, m ,\bar m,0,0) (\hat u, \hat v) = (z,g).
\end{equation*}
Moreover, there exists a suitable constant $\bar C$ such that
\begin{equation}
\label{invdiffestGD}
|\hat u| \le \bar C\Upsilon\Lambda| g|^{\upsilon'}_{\sigma, \mathbf{a}} + |z|^{\upsilon'}_{\sigma, \mathbf{b}}, \quad |\hat v| \le  | g|^{\upsilon'}_{\sigma, \mathbf{a}},
\end{equation}
where we recall that $|\hat u| = \max\{|\hat u|^{\upsilon'}_{\sigma, \mathbf{\bar b}}, | \left(\nabla \hat u \right) \Omega|^{\upsilon'}_{\sigma, \mathbf{b}} \}$ and $|\hat v| = \max\{|\hat v|^{\upsilon'}_{\sigma, \mathbf{\bar a}}, | \left(\nabla \hat v\right) \Omega|^{\upsilon'}_{\sigma, \mathbf{a}} \}$. Furthermore, $\Lambda$ is the constant in~\eqref{propabGD} and $\Upsilon$ is defined in~\eqref{H2GD}.
\end{lemma}
\begin{proof}
The proof of this lemma rests on Lemma \ref{homoeqlemmaGD}. Indeed, thanks to~\eqref{DFGD}, we can reformulate the problem in the following form. Given $(z,g) \in \mathcal{Z} \times \mathcal{G}$, we are looking for the unique solution $(\hat u, \hat v) \in \mathcal{U} \times \mathcal{V}$ of the following system
\begin{equation}
\label{systeminvGD}
\begin{cases}
\left(\nabla \hat u \right) \Omega = \bar m_0 \hat v - z\\
\left(\nabla \hat v \right) \Omega =g.
\end{cases} 
\end{equation}
By Lemma \ref{homoeqlemmaGD}, the unique solution $\hat v$ of the last equation of the latter system exists and satisfies
\begin{equation*}
|\hat v|^{\upsilon'}_{\sigma, \mathbf{\bar a}} \le |g |^{\upsilon'}_{\sigma, \mathbf{a}} . 
\end{equation*}
Moreover, by  $| \left(\nabla \hat v \right) \Omega|^{\upsilon'}_{\sigma, \mathbf{a}} = | g|^{\upsilon'}_{\sigma, \mathbf{a}}$, we have the second estimate in~\eqref{invdiffestGD}
\begin{equation}
\label{hatvGD}
|\hat v| = \max\{|\hat v|^{\upsilon'}_{\sigma, \mathbf{\bar a}}, | \left(\nabla \hat v \right) \Omega|^{\upsilon'}_{\sigma, \mathbf{a}} \} \le |g |^{\upsilon'}_{\sigma, \mathbf{a}}.
\end{equation}
Now, it remains to solve the first equation of~\eqref{systeminvGD} where $\hat v$ is known. For all fixed $t \in J_{\upsilon'}$ and thanks to property \textit{2.} of Proposition \ref{Holder}, the first condition of~\eqref{propabGD} and~\eqref{hatvGD}
\begin{eqnarray*}
{|(\bar m_0 \hat v - z)^t|_{C^\sigma} \over \mathbf{b}(t)} &\le& C\Upsilon {|\hat v^t|_{C^\sigma} \over \mathbf{b}(t)} + {|z^t|_{C^\sigma} \over \mathbf{b}(t)} \le C\Upsilon \Lambda{|\hat v^t|_{C^\sigma} \over \mathbf{\bar a}(t)} + {|z^t|_{C^\sigma} \over \mathbf{b}(t)} \\
&\le& C\Upsilon \Lambda |\hat v|^{\upsilon'}_{\sigma, \mathbf{ \bar a}} + |z|^{\upsilon'}_{\sigma, \mathbf{b}} \le C\Upsilon \Lambda |g |^{\upsilon'}_{\sigma, \mathbf{a}} + |z|^{\upsilon'}_{\sigma, \mathbf{b}},
\end{eqnarray*}
for a suitable constant $C$.
Taking the sup for all $t \in J_{\upsilon'}$ on the left-hand side of the latter, we obtain
\begin{equation*}
|\bar m_0 \hat v - z|^{\upsilon'}_{\sigma, \mathbf{b}}  \le C\Upsilon \Lambda |g |^{\upsilon'}_{\sigma, \mathbf{a}} + |z|^{\upsilon'}_{\sigma, \mathbf{b}}
\end{equation*}
and hence
\begin{equation*}
| \left(\nabla \hat u \right) \Omega|^{\upsilon'}_{\sigma, \mathbf{b}} = |\bar m_0 \hat v - z|^{\upsilon'}_{\sigma, \mathbf{b}}  \le  C\Upsilon \Lambda| g|^{\upsilon'}_{\sigma, \mathbf{a}} + |z|^{\upsilon'}_{\sigma, \mathbf{b}}.
\end{equation*}
Thanks to Lemma~\ref{homoeqlemmaGD} the unique solution $\hat u$ of the first equation of~\eqref{systeminvGD} exists verifying
\begin{equation*}
|\hat u|^{\upsilon'}_{\sigma, \mathbf{\bar b}} \le|\bar m_0 \hat v - z|^{\upsilon'}_{\sigma, \mathbf{b}(t)} \le C\Upsilon \Lambda|g |^{\upsilon'}_{\sigma, \mathbf{a}} + |z|^{\upsilon'}_{\sigma, \mathbf{b}}.
\end{equation*}
This concludes the proof of this lemma with $\bar C = C$ because
\begin{equation*}
|\hat u| = \max\{|\hat u|^{\upsilon'}_{\sigma, \mathbf{\bar b}}, | \left(\nabla \hat u \right) \Omega|^{\upsilon'}_{\sigma, \mathbf{b}} \} \le C\Upsilon \Lambda | g|^{\upsilon'}_{\sigma, \mathbf{a}} + |z|^{\upsilon'}_{\sigma, \mathbf{b}}.
\end{equation*}
\end{proof}

\subsection{$C^\sigma$-asymptotic KAM torus}\label{ICGD}
In the previous section, we proved that the functional $\mathcal{F}$ satisfies the hypotheses of the implicit function theorem. Here, we prove the existence of a $C^\sigma$-asymptotic KAM torus associated to $(X_H, X_{\tilde h}, \varphi_0)$ and we conclude the proof of  Theorem \ref{Thm1GD}.

Let $x = (a,b)$, where $a$ and $b$ are those defined by~\eqref{H2GD}. Obviously $(a,b) \in \mathcal{A} \times \mathcal{B}$ and 
\begin{equation}
\label{stimesimpaabGD}
|\partial_q a|^{\upsilon'}_{\sigma+1, \mathbf{a}} \le 1, \quad |b|^{\upsilon'}_{\sigma+2, \mathbf{b}} \le 1.
\end{equation}
We introduce the Banach space $(\mathcal{Y},| \cdot|)$ where $\mathcal{Y} = \mathcal{U} \times \mathcal{V}$ and, for all $y = (u,v) \in \mathcal{Y}$, $| y | = \max \{|u|, |v| \}$. Let $m$, $\bar m \in \mathcal{M}$ be as in~\eqref{H2GD} and we consider
\begin{equation*}
 \mathcal{F}(x,m, \bar m,y) = D_{(u,v)} \mathcal{F}(0,0,m, \bar m,0,0) y + \mathcal{R}(x,m, \bar m,y).
\end{equation*}
The aim is to find $y \in \mathcal{Y}$ in such a way that 
\begin{equation*}
\mathcal{F}(x,m, \bar m,y) =0,
\end{equation*}
where we recall that we have fixed $x$, $m$ and $\bar m$. This is equivalent to find $y \in \mathcal{Y}$ such that 
\begin{eqnarray*}
y &=&  -D_{(u,v)} \mathcal{F}(0,0,m, \bar m,0,0)^{-1}\mathcal{R}(x,m, \bar m,y)\\
&=& y - D_{(u,v)}  \mathcal{F}(0,0,m, \bar m,0,0)^{-1}\mathcal{F}(x,m, \bar m,y).
\end{eqnarray*}
This is well defined because we have already proved that $D_{(u,v)} \mathcal{F}(0,0,m, \bar m,0,0)$ is invertible (see Lemma \ref{lemmainvGD}). To this end, we introduce the following functional 
\begin{equation*}
\mathcal{L}(x,m, \bar m,\cdot) :  \mathcal{Y} \longrightarrow  \mathcal{Y}
\end{equation*}
 in such a way that 
\begin{equation}
\label{L}
\mathcal{L}(x,m, \bar m,y) =  y - D_{(u,v)} \mathcal{F}(0,0,m, \bar m,0,0)^{-1}\mathcal{F}(x,m, \bar m,y).
\tag{$\mathcal{L}$}
\end{equation}
This is well defined and, by the regularity of $\mathcal{F}$, we deduce that $\mathcal{L}$ is continuous, differentiable with respect to $y = (u,v)$ with differential $D_y\mathcal{L}$ continuous. The proof is reduced to find a fixed point of the latter. For this purpose, we introduce the following lemma. 

\begin{lemma}
\label{lemmautilethmGD}
There exists $\upsilon'$ large enough with respect to $n$, $\sigma$, $\Upsilon$, $\Lambda$ and $\mathbf{b}$, such that, for all  $y_*$,$y \in \mathcal{Y}$ with $|y_*|\le 1$,
\begin{equation}
\label{L2}
|D_y\mathcal{L}(x,m, \bar m,y_*) y| \le {1\over 2} |y|.
\end{equation}
\end{lemma}
\begin{proof}
The proof relies on Lemma \ref{lemmainvGD}. By~\eqref{L}, for all $y_*$,$y \in \mathcal{Y}$
\begin{align*}
&D_y\mathcal{L}(x,m, \bar m,y_*) y\\
&= D_{(u,v)} \mathcal{F}(0,0,m,\bar m,0,0)^{-1}\Big(D_{(u,v)}\mathcal{F}(0,0,m,\bar m,0,0) - D_{(u,v)}\mathcal{F}(x,m, \bar m, y_*)\Big)y.
\end{align*}
We can reformulate this problem in terms of estimating the unique solution $\hat y = (\hat u, \hat v) \in \mathcal{Y}$ of the following system
\begin{small}
\begin{equation}
\label{systultimolemma}
 D_{(u,v)} \mathcal{F}(0,0,m,\bar m,0,0)\hat y  = \Big(D_{(u,v)}\mathcal{F}(0,0,m,\bar m,0,0) - D_{(u,v)}\mathcal{F}(x,m,\bar m, y_*)\Big) y.
\end{equation}
\end{small}
Therefore, it suffices to estimate the right-hand side of the latter and apply Lemma \ref{lemmainvGD}. First, let us introduce the following notation. We observe that $y_* = (u_*, v_*) \in \mathcal{Y}$ and, for all $(q,t) \in \T^n \times J_{\upsilon'}$, we let
\begin{equation*}
\tilde u_*(q,t) = (q + u_*(q,t),t), \quad \tilde \varphi_*(q,t) = (q + u_*(q,t), v_*(q,t) ,t).
\end{equation*}
Thanks to~\eqref{DFGD}, the right-hand side of~\eqref{systultimolemma} is equal to 
\begin{small}
\begin{equation*}
\Big(D_{(u,v)}\mathcal{F}(0,0,m,\bar m,0,0) - D_{(u,v)}\mathcal{F}(x, m,\bar m, y_*)\Big) y = \begin{pmatrix} \bar m_0v - \left(\nabla u\right)\Omega - D_{(u,v)} F_1(b,\bar m, y_*)y\\
\left(\nabla v\right)\Omega - D_{(u,v)}F_2(x,m, y_*)y \end{pmatrix}
\end{equation*}
\end{small}
where, by~\eqref{DF1GD} and~\eqref{DF2GD},
\begin{eqnarray*}
 \bar m_0v - \left(\nabla u\right)\Omega - D_{(u,v)} F_1(b,\bar m, y_*)y &=& \left( \bar m_0 - \bar m \circ \tilde  \varphi_* \right) v - \left(\partial_q b \circ \tilde u_* \right)  u \\
&-& v_*^T \left(\partial_q \bar m \circ \tilde  \varphi_* \right) u - v_*^T \left( \partial_p \bar m \circ \tilde  \varphi_* \right) v \\
\left(\nabla v\right)\Omega - D_{(u,v)}F_2(x,m, y_*)y &=& -\left(\partial^2_q a \circ \tilde u_* \right)  u  - v_*^T \left( \partial^2_q b \circ \tilde u_* \right) u \\
&-&  (v_*^T)^2 \left( \partial^2_q m \circ \tilde  \varphi_* \right) u - \left(\partial_q b \circ \tilde u_* \right) v \\
&-&  (v_*^T)^2 \left(\partial^2_{pq} m \circ \tilde  \varphi_* \right)v - 2 v_*^T \left(\partial_q m \circ \tilde  \varphi_* \right)v.
\end{eqnarray*}
Thanks to property \textit{2.} of Proposition \ref{Holder}, we can estimate the first member on the left-hand side of the latter as follows
\begin{eqnarray*}
 \left|\left(\bar m_0v - \left(\nabla u\right)\Omega - D_{(u,v)} F_1(b,\bar m, y_*)y\right)^t\right|_{C^\sigma} &\le& C(\sigma)\Big( \left|\left(\bar m_0^t -\bar m \circ \tilde  \varphi_* \right)^t\right|_{C^\sigma} \left|v^t\right|_{C^\sigma}\\
&+& \left|\left(\partial_q b \circ \tilde u_* \right)^t\right|_{C^\sigma}  \left|u^t \right|_{C^\sigma}\\
&+& |v_*^t|_{C^\sigma} \left|\left(\partial_q \bar m \circ \tilde  \varphi_* \right)^t\right|_{C^\sigma} |u^t|_{C^\sigma} \\
&+&  |v_*^t|_{C^\sigma} \left|\left( \partial_p \bar m \circ \tilde  \varphi_* \right)^t\right|_{C^\sigma} |v^t|_{C^\sigma} \Big)
\end{eqnarray*}
for all $t \in J_{\upsilon'}$. 
We point out that $|y_*| = \max\{|u_*|, |v_*|\} \le 1$. Hence, we have to find an upper bound for each member on the right-hand side of the previous inequality. For all $t \in J_{\upsilon'}$
\begin{eqnarray*}
\left|\left(\bar m_0 -\bar m \circ \tilde  \varphi_* \right)^t\right|_{C^\sigma} |v^t|_{C^\sigma} &\le& C(\sigma) \Big(|\partial_q \bar m^t(\mathrm{id} + \tau u_*, \tau v_*)u^t_*|_{C^\sigma}\\
&+& |\partial_p \bar m^t(\mathrm{id} + \tau u_*, \tau v_*)v^t_*|_{C^\sigma}\Big) |v^t|_{C^\sigma} \\
&\le&C(\sigma)\Upsilon\left( 1 +  \mathbf{\bar b}(\upsilon')+  \mathbf{\bar a}(\upsilon')\right)|u^t_*|_{C^\sigma} |v^t|_{C^\sigma}\\
&+&C(\sigma)\Upsilon\left( 1 +  \mathbf{\bar b}(\upsilon')+  \mathbf{\bar a}(\upsilon')\right)|v^t_*|_{C^\sigma} |v^t|_{C^\sigma}\\
&\le&C(\sigma)\Upsilon\left(|u_*|\mathbf{\bar b}(t) +|v_*|\mathbf{\bar a}(t)\right)|v|\mathbf{\bar a}(t)\\
&\le&C(\sigma)\Upsilon\Lambda\mathbf{\bar b}(\upsilon')|y|\mathbf{b}(t) + C(\sigma)\Upsilon\Lambda^2\mathbf{b}(\upsilon')|y|\mathbf{b}(t)
\end{eqnarray*}
The first line of the latter is a consequence of the mean value theorem for a suitable $\tau \in [0,1]$. Concerning the second inequality, it is due to properties \textit{2.} and \textit{5.} of Proposition \ref{Holder}. Moreover, we use also that, thanks to~\eqref{propabGD} and for $\upsilon'$ large enough, we may assume $ \mathbf{\bar b}(\upsilon') \le 1$ and $ \mathbf{\bar a}(\upsilon') \le \Lambda \mathbf{b}(\upsilon') \le 1$. In the penultimate line on the right-hand side of the previous inequalities, we apply the following estimate  $ 1 +  \mathbf{\bar b}(\upsilon')+  \mathbf{\bar a}(\upsilon') \le 3$ for all $t \in J_{\upsilon'}$. In the last line, we use the first condition in~\eqref{propabGD}. 

Similarly to the previous case, thanks to property \textit{5.} of Proposition \ref{Holder}, the first condition in~\eqref{propabGD},~\eqref{stimesimpaabGD} and $ \mathbf{\bar b}(\upsilon') \le 1$, $ \mathbf{\bar a}(\upsilon') \le 1$, we obtain
\begin{eqnarray*}
\left|\left(\partial_q b \circ \tilde u_* \right)^t\right|_{C^\sigma}  \left|u^t \right|_{C^\sigma} &\le& C(\sigma) |b|^{\upsilon'}_{\sigma+2,\mathbf{b}} \mathbf{b}(t)\left( 1 +  \mathbf{\bar b}(\upsilon') \right)|u|\mathbf{\bar b}(t)\\
&\le&C(\sigma) \mathbf{\bar b}(\upsilon')|y| \mathbf{b}(t)\\
|v_*^t|_{C^\sigma} \left|\left(\partial_q \bar m \circ \tilde  \varphi_* \right)^t\right|_{C^\sigma} |u^t|_{C^\sigma} &\le&C(\sigma) |v_*|\mathbf{\bar a}(t) \Upsilon \left( 1 +  \mathbf{\bar b}(\upsilon')+  \mathbf{\bar a}(\upsilon')\right)|u|\mathbf{\bar b}(t)\\
&\le& C(\sigma)\Upsilon \Lambda\mathbf{\bar b}(\upsilon') |y| \mathbf{b}(t)\\
|v_*^t|_{C^\sigma} \left|\left( \partial_p \bar m \circ \tilde  \varphi_* \right)^t\right|_{C^\sigma} |v^t|_{C^\sigma} &\le&C(\sigma) |v_*|\mathbf{\bar a}(t) \Upsilon \left( 1 +  \mathbf{\bar b}(\upsilon')+  \mathbf{\bar a}(\upsilon')\right)|v|\mathbf{\bar a}(t)\\
&\le& C(\sigma)\Upsilon\Lambda^2 \mathbf{b}(\upsilon') |y| \mathbf{b}(t),
\end{eqnarray*}
for all $t \in J_{\upsilon'}$. Now, for $\upsilon'$ large enough, the previous estimates imply 
\begin{eqnarray*}
 \left|\left(\bar m_0v - \left(\nabla u\right)\Omega - D_{(u,v)} F_1(b,\bar m, y_*)y\right)^t\right|_{C^\sigma} &\le& {1 \over 4}|y|\mathbf{b}(t)
\end{eqnarray*}
for all $t \in J_{\upsilon'}$. Multiplying both sides of the latter by ${1 \over \mathbf{b}(t)}$ and taking the sup for all $t \in J_{\upsilon'}$, we obtain
\begin{eqnarray}
\label{FirstestimateGDfine}
 \left| \bar m_0v - \left(\nabla u\right)\Omega - D_{(u,v)} F_1(b,\bar m, y_*)y\right|^{\upsilon'}_{\sigma, \mathbf{b}} &\le&  {1 \over 4}|y|.
\end{eqnarray}
Similarly to the previous case, 
\begin{eqnarray*}
\left|\left(\left(\nabla v\right)\Omega - D_{(u,v)}F_2(x,m, y_*)y\right)^t\right|_{C^\sigma} &\le& C(\sigma)\Big(\left|\left(\partial^2_q a \circ \tilde u_* \right)^t\right|_{C^\sigma} \left| u^t\right|_{C^\sigma}\\
&+& \left|v_*^t\right|_{C^\sigma} \left|\left( \partial^2_q b \circ \tilde u_* \right)^t\right|_{C^\sigma} \left|u^t\right|_{C^\sigma} \\
&+&  \left|v_*^t\right|^2_{C^\sigma}\left| \left( \partial^2_q m \circ \tilde  \varphi_* \right)^t\right|_{C^\sigma} |u^t|_{C^\sigma} \\
&+&\left|\left(\partial_q b \circ \tilde u_* \right)^t\right|_{C^\sigma} \left|v^t\right|_{C^\sigma} \\
&+&  \left|v_*^t\right|_{C^\sigma}^2 \left|\left(\partial^2_{pq} m \circ \tilde  \varphi_* \right)^t\right|_{C^\sigma} \left|v^t\right|_{C^\sigma}\\ 
&+&  \left|v_*^t\right|_{C^\sigma} \left|\left(\partial_q m \circ \tilde  \varphi_* \right)^t\right|_{C^\sigma}\left|v^t\right|_{C^\sigma}\Big),
\end{eqnarray*}
for all $t \in J_{\upsilon'}$. Therefore, we have to estimate each member on the right-hand side of the latter. We begin with the element in the second line. For all $t \in J_{\upsilon'}$
\begin{eqnarray*}
|v_*^t|_{C^\sigma} \left|\left( \partial^2_q b \circ \tilde u_* \right)^t\right|_{C^\sigma} |u^t|_{C^\sigma} &\le& C(\sigma)|v_*^t|_{C^\sigma}|b|^{\upsilon'}_{\sigma+2,\mathbf{b}}\mathbf{b}(t)\left( 1 +  \mathbf{\bar b}(\upsilon') \right)|u^t|_{C^\sigma} \\
&\le& C(\sigma)\mathbf{\bar a}(t)|b|^{\upsilon'}_{\sigma+2,\mathbf{b}}\mathbf{b}(t) |u|\mathbf{\bar b}(t)\\
 &\le&  C(\sigma)\Lambda\mathbf{\bar b}(\upsilon')^2|y| \mathbf{a}(t).
\end{eqnarray*}
The first line of the above estimate is due to property \textit{5.} of Proposition \ref{Holder} and $ \mathbf{\bar b}(\upsilon') \le 1$. In the second line we use $|v_*^t|_{C^\sigma} \le \mathbf{\bar a}(t)$ and $|u^t|_{C^\sigma} \le |u|\mathbf{\bar b}(t)$ for all $t \in J_{\upsilon'}$. The last inequality is a consequence of the second condition of~\eqref{propabGD}. 

Thanks to property \textit{5.} of Proposition \ref{Holder},~\eqref{propabGD},~\eqref{stimesimpaabGD} and $ \mathbf{\bar b}(\upsilon') \le 1$, $ \mathbf{\bar a}(\upsilon') \le 1$, in the same way we have
\begin{eqnarray*}
\left|\left(\partial^2_q a \circ \tilde u_* \right)^t\right|_{C^\sigma} \left| u^t\right|_{C^\sigma} &\le& C(\sigma)|\partial_qa|^{\upsilon'}_{\sigma+1,\mathbf{a}} \mathbf{a}(t) |u|\mathbf{\bar b}(t) \\
&\le& C(\sigma)\mathbf{\bar b}(\upsilon') |y|\mathbf{a}(t),\\
\left|v_*^t\right|^2_{C^\sigma}\left| \left( \partial^2_q m \circ \tilde  \varphi_* \right)^t\right|_{C^\sigma} |u^t|_{C^\sigma} &\le& C(\sigma)|v_*|^2\mathbf{\bar a}(t)^2 \Upsilon |u|\mathbf{\bar b}(t) \\
&\le& C(\sigma)\Upsilon\Lambda \mathbf{\bar a}(t)\mathbf{b}(t) \mathbf{\bar b}(t)|u| \\
&\le&  C(\sigma)\Upsilon \Lambda^2\mathbf{\bar b}(\upsilon')^2 |y|\mathbf{a}(t)\\
\left|\left(\partial_q b \circ \tilde u_* \right)^t\right|_{C^\sigma} \left|v^t\right|_{C^\sigma} &\le& C(\sigma)|b|^{\upsilon'}_{\sigma+2,\mathbf{b}}\mathbf{b}(t)|v|\mathbf{\bar a}(t)\\
&\le& C(\sigma)\Lambda\mathbf{\bar b}(\upsilon') |y| \mathbf{a}(t)\\
 \left|v_*^t\right|_{C^\sigma}^2 \left|\left(\partial^2_{pq} m \circ \tilde  \varphi_* \right)^t\right|_{C^\sigma} \left|v^t\right|_{C^\sigma} &\le& C(\sigma) |v_*|^2\mathbf{\bar a}(t)^2 \Upsilon |v|\mathbf{\bar a}(t) \\
 &\le& C(\sigma) \Upsilon\Lambda^2\mathbf{b}(t)^2|v|\mathbf{\bar a}(t)\\
 &\le&  C(\sigma) \Upsilon\Lambda^3\mathbf{b}(\upsilon')\mathbf{\bar b}(\upsilon')|y|\mathbf{a}(t)\\
 \left|v_*^t\right|_{C^\sigma} \left|\left(\partial_q m \circ \tilde  \varphi_* \right)^t\right|_{C^\sigma}\left|v^t\right|_{C^\sigma} &\le& C(\sigma) |v_*|\mathbf{\bar a}(t) \Upsilon |v|\mathbf{\bar a}(t)\\
 &\le& C(\sigma)\Upsilon\Lambda\mathbf{b}(t)|v|\mathbf{\bar a}(t)\\
 &\le&C(\sigma)\Upsilon\Lambda^2\mathbf{\bar b}(\upsilon')|y|\mathbf{a}(t)\\
\end{eqnarray*}
for all $t \in J_{\upsilon'}$. Then, for $\upsilon'$ large enough
\begin{equation*}
\left|\left(\left(\nabla v\right)\Omega - D_{(u,v)}F_2(x,m, y_*)y\right)^t\right|_{C^\sigma} \le {1 \over 4\bar C\Upsilon \Lambda}|y|\mathbf{a}(t)
\end{equation*}
for all $t \in J_{\upsilon'}$. We recall that $\bar C$ is the constant introduced in Lemma \ref{lemmainvGD}. 
Multiplying both sides of the latter by ${1 \over \mathbf{a}(t)}$ and taking the sup for all $t \in J_{\upsilon'}$, we obtain
\begin{eqnarray}
\label{SecondestimateGDfine}
 \left|(\left(\nabla v\right)\Omega - D_{(u,v)}F_2(x,m, y_*)y\right|^{\upsilon'}_{\sigma, \mathbf{a}} &\le&  {1 \over 4\bar C\Upsilon\Lambda}|y|.
\end{eqnarray}
This concludes the proof of this lemma because, thanks to Lemma \ref{lemmainvGD}, the unique solution of~\eqref{systultimolemma} exists and by~\eqref{FirstestimateGDfine},~\eqref{SecondestimateGDfine}
\begin{eqnarray*}
|\hat u| &\le&  \left| \bar m_0v - \left(\nabla u\right)\Omega - D_{(u,v)} F_1(b,\bar m, y_*)y\right|^{\upsilon'}_{\sigma, \mathbf{b}} \\
&+&  \bar C\Upsilon\Lambda \left|(\left(\nabla v\right)\Omega - D_{(u,v)}F_2(x,m, y_*)y\right|^{\upsilon'}_{\sigma, \mathbf{a}}  \le {1 \over 2}|y|\\
|\hat v| &\le&  \left|(\left(\nabla v\right)\Omega - D_{(u,v)}F_2(x,m, y_*)y\right|^{\upsilon'}_{\sigma, \mathbf{a}}  \le {1 \over 4\bar C\Upsilon \Lambda}|y| \le {1\over 2}|y|.
\end{eqnarray*}
\end{proof}

We observe that the choice of the constant $1$ in the ball $|y_*|\le 1$ is completely arbitrary. One can choose another threshold provided to take $\upsilon'$ sufficiently large. 

Now, the previous lemma proves that $\mathcal{L}(x,m, \bar m,\cdot)$ is a contraction of a complete subset of $\mathcal{Y}$. Then, there exists a unique fixed point $y \in \mathcal{Y}$ with $|y| \le 1$. This concludes the proof of Theorem \ref{Thm1GD}. 

\section{Proof of Corollary \ref{Cor1GD}}\label{ProofC1GD}
The proof is essentially the same as that of Theorem \ref{Thm1GD}. Because of that, we will not give all the details. However, we will provide the necessary elements to reconstruct the proof.

We are looking for a $C^\sigma$-asymptotic KAM torus $\psi^t$ associated to $(Z, \omega, \mathrm{Id})$, where $Z$ is the vector field defined by~\eqref{ZGDIntro}. This means that, for given $Z$, we are searching for $\upsilon' \ge 0$ sufficiently large and a suitable function $u : \T^n \times J_{\upsilon'} \to \R^n$ such that 
\begin{equation*}
\psi(q,t) = q + u(q,t)
\end{equation*}
and, in addition, $\psi$ and $u$ satisfy 
\begin{align}
\label{hyp1GD2Z}
&  Z(\psi(q, t), t) -  \partial_q \psi(q, t) \omega - \partial_t \psi(q, t) = 0,\\
\label{hyp2GD2Z}
&  \lim_{t \to +\infty}  |u^t|_{C^\sigma} = 0
\end{align}
for all $(q,t) \in \T^n \times J_{\upsilon'}$. We will choose $\upsilon'$ sufficiently large in Lemma \ref{lemmaCorA}.
Similarly to the proof of Theorem \ref{Thm1GD}, we introduce a suitable functional $\mathcal{F}$ given by~\eqref{hyp1GD2Z}. To this end, we define
\begin{equation*}
\tilde \psi(q,t) = (q + u(q,t), t), 
\end{equation*}
for all $(q,t) \in \T^n \times J_{\upsilon'}$. The composition of $Z$ with $\tilde \psi$ is equal to
\begin{equation*}
Z \circ \tilde \psi (q, t) =\omega + P \circ \tilde \psi (q, t) 
\end{equation*}
and
\begin{equation*}
\partial_q \psi(q, t)\omega + \partial_t \psi(q, t) = \omega + \partial_q u(q, t)\omega + \partial_t u (q, t)
\end{equation*}
for all $(q, t) \in \T^n \times J_{\upsilon'}$.
We recall the notation introduced in the previous section 
\begin{equation*}
\nabla u (q, t) \Omega = \partial_q u(q, t) \omega + \partial_t u(q, t)
\end{equation*}
for all $(q, t) \in \T^n \times J_{\upsilon'}$.
Then,  we can rewrite~\eqref{hyp1GD2Z} in the following form 
\begin{eqnarray}
\label{FGDZC}
P\circ \tilde\psi - \left(\nabla u \right) \Omega = 0,
\end{eqnarray}
where we have omitted the arguments $(q, t)$.  

Before the introduction of the functional $\mathcal{F}$, let $\upsilon \ge 0$ and $\sigma \ge 1$ be the positive parameters defined in Corollary \ref{Cor1GD}. For $\upsilon' \ge \upsilon \ge0$ that will be chosen later, we introduce the following Banach spaces $\left(\mathcal{P}, |\cdot |\right)$, $\left(\mathcal{U}, |\cdot |\right)$ and $\left(\mathcal{Z}, |\cdot |\right)$  
\vspace{5mm}
\begin{eqnarray*}
\mathcal{P} &=& \Big\{P : \T^n \times J_{\upsilon'}  \to \R^n  \hspace{1mm}| \hspace{1mm} P \in \mathcal{\bar S}^{\upsilon'}_{\sigma, 1} , \hspace{1mm} \mbox{and} \hspace{1mm} |P| =|P|^{\upsilon'}_{\sigma+1, \mathbf{P}} < \infty\Big\}\\
\mathcal{U} &=& \Big\{u :\T^n \times J_{\upsilon'}  \to \R^n  \hspace{1mm}| \hspace{1mm} u,\left( \nabla u\right) \Omega \in \mathcal{S}^{\upsilon'}_{\sigma} \\
&& \mbox{and} \hspace{1mm} |u| = \max\{|u|^{\upsilon'}_{\sigma, \mathbf{\bar P}}, | \left(\nabla u \right) \Omega|^{\upsilon'}_{\sigma, \mathbf{P}} \} < \infty\Big\}\\
\mathcal{Z} &=& \Big\{z : \T^n \times J_{\upsilon'}  \to \R^n  \hspace{1mm}| \hspace{1mm} z \in \mathcal{S}^{\upsilon'}_{\sigma}   , \hspace{1mm} \mbox{and} \hspace{1mm} |z| =|z|^{\upsilon'}_{\sigma, \mathbf{P}} < \infty\Big\}
\end{eqnarray*} 
Let $\mathcal{F}$ be the following functional
\begin{equation*}
\mathcal{F} :  \mathcal{P} \times  \mathcal{U} \longrightarrow \mathcal{Z} \\
\end{equation*}
\begin{equation*}
\mathcal{F}(P,u) = P\circ \tilde\psi - \left(\nabla u \right) \Omega.
\end{equation*}
This is obtained by~\eqref{FGDZC} and we observe that 
\begin{equation*}
\mathcal{F}(0, 0)=0.
\end{equation*}
We can reformulate our problem in the following form. For $P \in \mathcal{P}$ sufficiently close to $0$, we are looking for $u \in \mathcal{U}$ in such a way that $\mathcal{F}(P,u) = 0$.  

Concerning the associated linearized problem,  the differential of $\mathcal{F}$ with respect to the variable $u$ calculated in $(0,0)$ is equal to
\begin{equation*}
D_u \mathcal{F}(0,0)\hat u = -\left(\nabla \hat u \right) \Omega.
\end{equation*}

The functional $\mathcal{F}$ is well defined, continuous, differentiable with respect to $u$ with $D_u\mathcal{F}(P,u)$ continuous. Moreover, as a straightforward consequence of Lemma \ref{homoeqlemmaGD}, $D_u \mathcal{F}(0,0)$ is invertible. Then, $\mathcal{F}$ satisfies the hypotheses of the implicit function theorem. 

Now, similarly to the proof of Theorem \ref{Thm1GD}, we fix $P$ as in Corollary \ref{Cor1GD} and we introduce the following functional
\begin{equation*}
\mathcal{L}(P,\cdot):\mathcal{U} \longrightarrow \mathcal{U}
\end{equation*}
 in such a way that 
\begin{equation*}
\mathcal{L}(P,u) =  u - D_u \mathcal{F}(0,0)^{-1}\mathcal{F}(P,u).
\end{equation*}
We recall that $P$ is fixed and the proof of Corollary \ref{Cor1GD} is reduced to find a fixed point of the latter. To this end, we have the following lemma
\begin{lemma}
\label{lemmaCorA}
There exists  $\upsilon'$ large enough with respect to $n$, $\sigma$ and $\mathbf{P}$, such that, for all  $u_*$,$u \in \mathcal{U}$ with $|u_*|\le 1$,
\begin{equation*}
|D_u\mathcal{L}(P,u_*) u| \le {1\over 2} |u|.
\end{equation*}
\end{lemma}

Therefore, $\mathcal{L}(P, \cdot)$ is a contraction of a complete subset of $\mathcal{P}$ and this concludes the proof of Corollary \ref{Cor1GD}.

\section{Proof of Theorem \ref{Thm1AGD}}\label{ProofT1AGD}

The proof of this theorem is extremely similar to that of Theorem \ref{Thm1GD}. 
We expand the Hamiltonian $H$ in~\eqref{H1AGD} in a small neighbourhood of $0 \in B$. Then, thanks to Proposition \ref{propAAGD} and for a positive parameter $\Upsilon \ge 1$, we can rewrite the Hamiltonian $H$ in the following form
\begin{equation}
\label{H2AGD}
\begin{cases}
H : \T^n \times B \times J_0 \longrightarrow \R\\
H(q,p,t) = \omega \cdot p + a(q,t) + b(q,t) \cdot p + m(q,p,t) \cdot p^2,\\
a, b, \partial^2_p H \in \mathcal{A}^0_s,\\
\sup_{t \in J_0} |a^t|_s < \infty, \quad \sup_{t \in J_0}|\partial_p^2 H^t|_s \le \Upsilon,\\
|\partial_q a^t|_s \le \mathbf{a}(t), \quad |b^t|_s  \le \mathbf{b}(t), \quad \mbox{for all $t \in J_0$}
\end{cases}
\tag{$**_B$}
\end{equation}
where $s = {s_0 \over 2}$ and $|\cdot |_s$ is the analytic norm. Furthermore, $\mathbf{a}(t)$ and $\mathbf{b}(t)$ are the functions introduced in~\eqref{H1AGD} satisfying~\eqref{propabGD} and $\mathcal{A}^0_s$ is the space of functions defined in Definition \ref{AsetfunctionGDA}.

We consider  the following Hamiltonian 
\begin{equation*}
\tilde h(q,p,t) = h(q,p,t) + \int_0^1 (1 -\tau)\partial^2_p f(q, \tau p, t) d\tau\cdot p^2.
\end{equation*}
for all $(q,p,t) \in \T^n \times B \times J_0$. It is obvious that $\tilde h \in \mathcal{K}_\omega$ and $X_H$, $X_{\tilde h}$ verify~\eqref{hyp1introAKAM}.

\subsection{Preliminary Settings}\label{PSAGD}

We are looking for an analytic asymptotic KAM torus $\varphi^t$ associated to $(X_H, X_{\tilde h}, \varphi_0)$, where $H$ is the Hamiltonian in~\eqref{H2AGD}, $\tilde h$ is the Hamiltonian previously defined and $\varphi_0$ the trivial embedding $\varphi_0 : \T^n  \to \T^n \times B$, $\varphi_0(q) = (q,0)$. More specifically, for given $H$, we are searching for $\upsilon' \ge 0$ sufficiently large and suitable functions $u$, $v : \T^n \times J_{\upsilon'} \to \R^n$ such that 
\begin{equation*}
\varphi(q,t) = (q + u(q,t), v(q,t))
\end{equation*}
and in such a way that $\varphi$, $u$ and $v$ satisfy
\begin{align}
\label{hyp1AGD2}
&   X_H(\varphi(q, t), t) -  \partial_q \varphi(q, t) \omega - \partial_t \varphi(q, t) = 0,\\
\label{hyp2AGD2}
&  \lim_{t \to +\infty}  |u^t|_{s \over 2}= 0, \quad \lim_{t \to +\infty}  |v^t|_{s \over 2} = 0,
\end{align}
for all $(q,t) \in \T^n \times J_{\upsilon'}$. Similarly to Theorem \ref{Thm1GD}, we will choose $\upsilon'$  large enough in Lemma \ref{lemmautilethmAGD} ($\upsilon'$  will be already required large in Lemma \ref{lemmaquantAGD}).

To this end, we introduce a special norm and suitable Banach spaces. Given $s>0$ and  $\upsilon \ge 0$, for every $f \in \mathcal{A}^\upsilon_s$ and for positive real functions $\mathbf{u}(t)$ defined on $J_\upsilon$, we consider the following norm
\begin{equation*}
|f|^\upsilon_{s, \mathbf{u}} = \sup_{t \in J_\upsilon}{|f^t|_s \over \mathbf{u}(t)}.
\end{equation*}
It is the analytic version of the norm defined in the finitely differentiable case by~\eqref{defnormGD}. 

Let $\upsilon \ge 0$, $s>0$ and $\Upsilon \ge 1$ be the positive parameters introduced by~\eqref{H2AGD} and~\eqref{propabGD}. For $\upsilon' \ge \upsilon \ge 0$ that will be chosen later, we consider the following Banach spaces $\left(\mathcal{A}, |\cdot |\right)$, $\left(\mathcal{B}, |\cdot |\right)$, $\left(\mathcal{U}, |\cdot |\right)$, $\left(\mathcal{V}, |\cdot |\right)$, $\left(\mathcal{Z}, |\cdot |\right)$, $\left(\mathcal{G}, |\cdot |\right)$ and  $\left(\mathcal{M}, |\cdot|\right)$. 
\vspace{5mm}
\begin{eqnarray*}
\mathcal{A} &=& \Big\{a : \T^n \times J_{\upsilon'} \to \R  \hspace{1mm}| \hspace{1mm} a \in \mathcal{A}^{\upsilon'}_s \hspace{1mm} \mbox{and} \hspace{1mm} |a| =|a|^{\upsilon'}_{s, 1} + |\partial_q a|^{\upsilon'}_{s, \mathbf{a}} < \infty\Big\}\\
\mathcal{B} &=& \Big\{b : \T^n \times J_{\upsilon'}  \to \R^n  \hspace{1mm}| \hspace{1mm} b \in \mathcal{A}^{\upsilon'}_s , \hspace{1mm} \mbox{and} \hspace{1mm} |b| =|b|^{\upsilon'}_{s, \mathbf{b}} < \infty\Big\}\\
\mathcal{U} &=& \Big\{u :\T^n \times J_{\upsilon'}  \to \R^n  \hspace{1mm}| \hspace{1mm} u,\left( \nabla u\right) \Omega \in \mathcal{A}^{\upsilon'}_{s \over 2} \\
&& \mbox{and} \hspace{1mm} |u| = \max\{|u|^{\upsilon'}_{{s \over 2}, \mathbf{\bar b}}, | \left(\nabla u \right) \Omega|^{\upsilon'}_{{s \over 2}, \mathbf{b}} \} < \infty\Big\}\\
\mathcal{V} &=& \Big\{v :\T^n \times J_{\upsilon'}  \to \R^n  \hspace{1mm}| \hspace{1mm} v,\left( \nabla v\right) \Omega \in \mathcal{A}^{\upsilon'}_{s \over 2} \\
&& \mbox{and} \hspace{1mm} |v| = \max\{|v|^{\upsilon'}_{{s \over 2}, \mathbf{\bar a}}, |\left( \nabla v \right) \Omega|^{\upsilon'}_{{s \over 2}, \mathbf{a}} \} < \infty\Big\}\\
\mathcal{Z} &=& \Big\{z : \T^n \times J_{\upsilon'}  \to \R^n  \hspace{1mm}| \hspace{1mm} z \in \mathcal{A}^{ \upsilon'}_{s \over 2} , \hspace{1mm} \mbox{and} \hspace{1mm} |z| =|z|^{\upsilon'}_{{s \over 2}, \mathbf{b}} < \infty\Big\}\\\mathcal{G} &=& \Big\{g : \T^n \times J_{\upsilon'} \to \R  \hspace{1mm}| \hspace{1mm} g \in \mathcal{A}^{ \upsilon'}_{s \over 2} \hspace{1mm} \mbox{and} \hspace{1mm} |g| =|g|^{\upsilon'}_{{s \over 2}, \mathbf{a}} < \infty\Big\}\\
\mathcal{M} &=& \Big\{m : \T^n \times B \times J_{\upsilon'} \to M_n  \hspace{1mm}| \hspace{1mm} m \in \mathcal{A}^{\upsilon'}_s \hspace{1mm} \mbox{and} \hspace{1mm} |m| =|m|^{\upsilon'}_{s, 1} \le \Upsilon\Big\}
\end{eqnarray*} 
The risk of mixing the Banach space $\mathcal{A}$ with the space of functions $\mathcal{A}^\upsilon_s$ is small.
Concerning $\mathcal{A}$, we have that $|a|^{\upsilon'}_{s, 1} = \sup_{t \in J_{\upsilon'}}|a^t|_s$, similarly for $\mathcal{M}$. Regarding the last Banach space $\mathcal{M}$, we recall that $M_n$ is the set of $n$-dimensional matrices. These are the analytic version of the Banach spaces introduced in the finitely differentiable case (see  Section \ref{PrelSettGD}).

Now, let us introduce the following subspace $\mathcal{X}$ of $\mathcal{A} \times \mathcal{B} \times \mathcal{M} \times \mathcal{M} \times  \mathcal{U} \times  \mathcal{V}$, in such a way that 
\begin{eqnarray*}
\mathcal{X} &=& \{(a,b,m, \bar m,u,v) \in \mathcal{A} \times \mathcal{B} \times \mathcal{M} \times \mathcal{M} \times  \mathcal{U} \times  \mathcal{V} :\\
&& |\partial_q a|^{\upsilon'}_{s, \mathbf{a}} \le 1, \hspace{2mm}  |b| \le 1, \hspace{2mm} |m| \le \Upsilon, \hspace{2mm} |\bar m| \le \Upsilon, \hspace{2mm}  |u| \le 1, \hspace{2mm}  |v| \le  1\}.
\end{eqnarray*}

Let $\mathcal{F}$ be the following functional
\begin{equation}
\label{FfunctionalAGD}
\mathcal{F} :  \mathcal{X} \longrightarrow \mathcal{Z} \times \mathcal{G}\\
\end{equation}
\begin{equation*}
\mathcal{F}(a, b,m , \bar m,  u, v) = (F_1(b, \bar m, u,v), F_2(a,b,m,u,v))
\end{equation*}
with 
\begin{eqnarray*}
F_1(b,\bar m ,u,v) &=& b \circ \tilde u + \left(\bar m \circ \tilde \varphi \right) v - \left(\nabla u \right)\Omega,\\
F_2(a,b, m,u,v) &=& \partial_q a\circ \tilde u  + \left(\partial_q b \circ \tilde u \right) v  + \left(\partial_q m \circ \tilde \varphi \right) \cdot v^2 + \left(\nabla v \right)\Omega.
\end{eqnarray*}

Similarly to Section \ref{OPGD}, this functional is obtained by~\eqref{hyp1AGD2}. Therefore, contrary to the proof of Theorem \ref{Thm1GD}, we have to define $\mathcal{F}$ on a suitable subspace $\mathcal{X}$ of $\mathcal{A} \times \mathcal{B} \times \mathcal{M} \times \mathcal{M} \times  \mathcal{U} \times  \mathcal{V}$. This is because we have to control the domain of analyticity of the components of $\mathcal{F}$. We have to verify that $\mathcal{F}$ satisfies the hypotheses of the implicit function theorem. 

\subsection{Homological Equation}

Given $s > 0$, $\upsilon \ge 0$ and $\omega \in \R^n$, we are looking for a solution of  the following equation for the unknown $\varkappa : \T_s^n \times J_\upsilon \to \R$
\begin{equation}
\label{HEAGD}
\begin{cases}
 \omega \cdot \partial_q \varkappa(q,t) + \partial_t \varkappa(q,t) = g(q,t),\\
g \in \mathcal{A}^\upsilon_s, \quad |g|_{s, \mathbf{g}}^\upsilon < \infty
\end{cases}
\tag{$HE_B$}
\end{equation}
where $\mathbf{g}(t)$ is a positive, decreasing, integrable function on $J_\upsilon$ and $g : \T_s^n \times J_\upsilon \to \R$ is given. 

\begin{lemma}[\textbf{Homological Equation}]
\label{homoeqlemmaAGD} 
There exists a unique solution $\varkappa \in  \mathcal{A}^\upsilon_s$ of~\eqref{HEAGD} such that 
\begin{equation*}
\lim_{t \to +\infty} |\varkappa^t|_{C^0} = 0.
\end{equation*}
Moreover, 
\begin{equation*}
|\varkappa|^\upsilon_{s, \mathbf{\bar g}} \le |g|_{s, \mathbf{g}}^\upsilon.
\end{equation*}
\end{lemma}
\begin{proof}
The proof of this lemma is essentially the same as that of Lemma \ref{homoeqlemmaGD}.  

\textit{Existence}: Let us define the following transformation 
\begin{equation*}
h(q,t) = (q -  \omega t , t),
\end{equation*}
then $h: \T^n_s \times J_\upsilon \to \T^n_s \times J_\upsilon$ because $ \omega \in \R^n$, $t \in J_\upsilon \subset \R$ and thus $q -  \omega t \in \T_s^n$ if and only if $q \in \T^n_s$. The fact that $t$ is real and not complex is of fundamental importance. This ensures that the latter transformation is well-defined. 

It is enough to prove the first part of this lemma for the much simpler equation
\begin{equation*}
\partial_t \kappa = g(q+ \omega t, t).
\end{equation*}
The unique solution of the above equation satisfying the asymptotic condition is 
\begin{equation*}
\kappa(q,t) = -\int_t^{+\infty} g(q+ \omega\tau, \tau)d\tau
\end{equation*}
and hence
\begin{equation*}
\varkappa(q,t) = \kappa \circ h(q,t) =-\int_t^{+\infty} g(q+ \omega(\tau - t), \tau)d\tau
\end{equation*}
is the unique solution of~\eqref{HEAGD} that we are looking for.

\textit{Regularity and Estimates}: $g \in \mathcal{A}^\upsilon_s$ implies $\kappa \in \mathcal{A}^\upsilon_s$ and hence $\varkappa = \kappa \circ h \in\mathcal{A}^\upsilon_s$. Moreover, for all fixed $t \in J_\upsilon$
\begin{equation*}
|\varkappa^t|_s \le |g|_{s, \mathbf{g}}^\upsilon \mathbf{\bar g}(t).
\end{equation*}
We prove the second part of this lemma by multiplying both sides of the latter by ${1 \over \mathbf{\bar g}(t)}$ and taking the sup for all $t \in J_\upsilon$.
\end{proof}

\subsection{Regularity of $\mathcal{F}$}

We begin this part with the following quantitative lemma. It is of fundamental importance to prove that $\mathcal{F}$ (see~\eqref{FfunctionalAGD}) is well defined. The following lemma imposes the first restriction on $\upsilon'$. We will take a stronger one after.

\begin{lemma}
\label{lemmaquantAGD}
For $\upsilon'$ large enough with respect to $s$, $\Lambda$ and $\mathbf{b}$, if $(u,v) \in \mathcal{U} \times \mathcal{V}$ satisfies the following estimates $|u| \le 1$ and $|v| \le 1$, then
\begin{equation*}
\sup_{t \in J_{\upsilon'}}|u^t|_{s\over 2} \le {s \over 8}, \quad \sup_{t \in J_{\upsilon'}}|v^t|_{s\over 2} \le {s \over 8}.
\end{equation*}
\end{lemma}
\begin{proof}
 If $|u| \le 1$ and $|v| \le 1$, then by~\eqref{propabGD}
\begin{equation*}
|u^t|_{s\over 2} \le \mathbf{\bar b}(t) \le \mathbf{\bar b}(\upsilon'), \quad |v^t|_{s\over 2} \le \mathbf{\bar a}(t) \le \Lambda \mathbf{b}(t)  \le \Lambda \mathbf{b}(\upsilon').
\end{equation*}
for all $t \in J_{\upsilon'}$. Now, for $\upsilon'$ large enough, we have the claim.
\end{proof}

Thanks to the previous lemma, the properties in Appendix \ref{B} and~\eqref{propabGD}, one can prove that $\mathcal{F}$ is well defined. Moreover, similarly to Theorem \ref{Thm1GD}, it is continuous, differentiable with respect to the variables $(u,v)$ and $D_{(u,v)}\mathcal{F}$ is continuous. We recall that $D_{(u,v)}$ stands for the differential with respect to $(u,v)$. Furthermore, $D_{(u,v)}\mathcal{F}$ calculated on $(0,0,m, \bar m, 0, 0)$ is equal to 
\begin{equation*}
D_{(u,v)} \mathcal{F}(0,0,m,\bar m, 0, 0)(\hat u, \hat v) = (\bar m_0\hat v - \left(\nabla \hat u\right)\Omega,  \left(\nabla \hat v\right)\Omega),
\end{equation*}
where, for all $(q, t) \in \T^n \times J_{\upsilon'}$, $\bar m_0(q,t) = \bar m(q,0,t)$. As one can expect, for all fixed $m$, $\bar m \in \mathcal{M}$, $D_{(u,v)}\mathcal{F}(0,0, m, \bar m,0,0) $ is invertible. More specifically, we have the following lemma

\begin{lemma}
\label{lemmainvAGD}
For all $(z,g) \in \mathcal{Z} \times \mathcal{G}$ there exists a unique $(\hat u, \hat v) \in \mathcal{U} \times \mathcal{V}$ such that 
\begin{equation*}
D_{(u,v)}\mathcal{F}(0,0,m,\bar m,0,0) (\hat u, \hat v) = (z,g).
\end{equation*}
Moreover, for a suitable constant $\bar C$
\begin{equation*}
|\hat u| \le \bar C \Upsilon \Lambda | g|^{\upsilon'}_{{s\over 2}, \mathbf{a}} + |z|^{\upsilon'}_{{s\over 2}, \mathbf{b}}, \quad |\hat v| \le  | g|^{\upsilon'}_{{s\over 2}, \mathbf{a}}.
\end{equation*}
\end{lemma}
\begin{proof}
The proof is essentially the same as that of Lemma \ref{lemmainvGD}. It relies on Lemma \ref{homoeqlemmaAGD}.
\end{proof}

\subsection{Analytic asymptotic KAM torus}

Let $a$ and $b$ be the functions introduced by~\eqref{H2AGD}. It is straightforward to verify that $(a,b)\in\mathcal{A} \times \mathcal{B}$ and
\begin{equation}
\label{stimesimpaabAGD}
|\partial_q a|^{\upsilon'}_{s, \mathbf{a}} \le 1, \quad |b|^{\upsilon'}_{s, \mathbf{b}} \le 1.
\end{equation}
We introduce the Banach space $(\mathcal{Y},| \cdot|)$, such that $\mathcal{Y} = \mathcal{U} \times \mathcal{V}$ and for all $y = (u,v) \in \mathcal{Y}$, $| y | = \max \{|u|, |v| \}$. Following the lines of the differentiable case (Section~\eqref{ICGD}), we fix $m$, $\bar m \in \mathcal{M}$ as in~\eqref{H2AGD}, we let $x = (a,b)$ and we introduce the following functional
\begin{equation*}
\mathcal{L}(x,m, \bar m,\cdot) :  \mathcal{Y} \longrightarrow  \mathcal{Y}
\end{equation*}
 in such a way that 
\begin{equation}
\label{LAGD}
\mathcal{L}(x,m, \bar m,y) =  y - D_{(u,v)} \mathcal{F}(0,0,m, \bar m,0,0)^{-1}\mathcal{F}(x,m, \bar m,y).
\tag{$\mathcal{L}$}
\end{equation}
This is well defined and, by the regularity of $\mathcal{F}$, we deduce that $\mathcal{L}$ is continuous, differentiable with respect to $y = (u,v)$ with differential $D_y\mathcal{L}$ continuous. The proof is reduced to find a fixed point of the latter. For this purpose, we state the following lemma, which is the analytic version of Lemma \ref{lemmautilethmGD}. 

\begin{lemma}
\label{lemmautilethmAGD}
There exists $\upsilon'$ large enough with respect to $s$, $\Upsilon$, $\Lambda$ and $\mathbf{b}$, such that, for all  $y_*$,$y \in \mathcal{Y}$ with $|y_*|\le 1$,
\begin{equation*}
|D_y\mathcal{L}(x,m, \bar m,y_*) y| \le {1\over 2} |y|.
\end{equation*}
\end{lemma}
\begin{proof}
The proof of this lemma is extremely similar to that of Lemma \ref{lemmautilethmGD}. For this reason, it is omitted.
\end{proof}

The above lemma proves that $\mathcal{L}(x, m, \bar m, \cdot)$ is a contraction and this concludes the proof of Theorem \ref{Thm1AGD}.


\appendix
\section{Hölder classes of functions}
\label{A}
This part is dedicated to a very brief introduction to Hölder classes of functions $C^\sigma$. Let $D$ be an open subset of $\R^n$. For integers $k \ge 0$, we denote by $C^k(D)$ the spaces of functions $f : D \to \R$ with continuous partial derivatives $\partial^\alpha f \in C^0(D)$ for all $\alpha \in \N^n$ with $|\alpha|=\alpha_1+...+\alpha_n \le k$. We define the norm
\begin{equation*}
|f|_{C^k} = \sup_{|\alpha|\le k}|\partial^\alpha f|_{C^0},
\end{equation*}
where $|\partial^\alpha f|_{C^0} = \sup_{x \in D}|\partial^\alpha f(x)|$ denotes the sup norm. For $\sigma=k+\mu$, with $k \in \Z$, $k \ge 0$ and $0 < \mu <1$, the Hölder spaces $C^\sigma(D)$ are the spaces of functions $f\in C^k(D)$ such that $|f|_{C^\sigma}  <\infty$, where
\begin{equation}
\label{Holdernorm}
|f|_{C^\sigma} = \sup_{|\alpha|\le k}|\partial^\alpha f|_{C^0} + \sup_{|\alpha| = k}{|\partial^\alpha f(x) - \partial^\alpha f(y)| \over |x-y|^\mu}.
\end{equation}
In the case of functions $f=(f_1,...,f_n)$ with values in $\R^n$, we set $|f|_{C^\sigma} = \max_{1 \le i \le n} |f_i|_{C^\sigma}$. Moreover, in agreement with the convention made above, if $M =\{m_{ij}\}_{1 \le i,j \le n}$ is a $n \times n$ matrix, we set $|M|_{C^\sigma}= \max_{1 \le i,j \le n} |m_{ij}|_{C^\sigma}$. 

In what follows, we have some properties of these norms that we widely use in this work. First, we recall that $C(\cdot)$ stands for constants depending on  $n$ and other parameters into brackets.
\begin{proposition}
\label{convexity}
For all $f\in C^{\sigma_1}(\R^n)$, then 
\begin{equation*}
|f|^{\sigma_1-\sigma_0}_{C^\sigma}\le C(\sigma_1) |f|^{\sigma_1-\sigma}_{C^{\sigma_0}}|f|^{\sigma-\sigma_0}_{C^{\sigma_1}} \hspace{3mm} \mbox{for all $0 \le \sigma_0 \le \sigma \le \sigma_1$}.
\end{equation*}
\end{proposition}
\begin{proof}
We refer to~\cite{Hor76} for the proof. 
\end{proof}
Furthermore, we have the following Proposition. 
\begin{proposition}
\label{Holder}
We consider $f$, $g \in C^\sigma(D)$ and $\sigma \ge 0$.
\begin{enumerate}
\item For all $\beta \in \N^{n}$, if $|\beta| + s = \sigma$ then  $\left|{\partial^{|\beta|} \over \partial{x_1}^{\beta_1}... \partial{x_n}^{\beta_n}} f \right|_{C^s} \le |f|_{C^\sigma}$.\\
\item  $|fg|_{C^\sigma} \le C(\sigma)\left(|f|_{C^0}|g|_{C^\sigma} + |f|_{C^\sigma}|g|_{C^0}\right)$. 
\end{enumerate}
Now we consider composite functions. Let $z$ be defined on $D_1 \subset \R^n$ and takes its values on $D_2 \subset \R^n$ where $f$ is defined. 

If $\sigma < 1$, $f \in C^1(D_2)$, $z \in C^\sigma (D_1)$ then $f\circ z \in C^\sigma(D_1)$ 
\begin{enumerate}
\item[3.] $|f \circ z|_{C^\sigma} \le C(|f|_{C^1}|z|_{C^\sigma}+ |f|_{C^0})$.
\end{enumerate}
If $\sigma < 1$, $f \in C^\sigma(D_2)$, $z \in C^1 (D_1)$ then $f\circ z \in C^\sigma(D_1)$  
\begin{enumerate}
\item[4.] $|f \circ z|_{C^\sigma} \le C(|f|_{C^\sigma}|\nabla z|^\sigma_{C^0}+ |f|_{C^0})$.
\end{enumerate}
If $\sigma \ge 1$ and $f \in C^\sigma (D_2)$, $z \in C^\sigma (D_1)$ then $f\circ z \in C^\sigma(D_1)$ 
\begin{enumerate}
\item[5.] $|f \circ z|_{C^\sigma} \le C(\sigma) \left(|f|_{C^\sigma}|\nabla z|^\sigma_{C^0} + |f|_{C^1}|\nabla z|_{C^{\sigma-1}}+ |f|_{C^0}\right)$.
\end{enumerate}
\end{proposition}
\begin{proof}
The proofs of the properties contained in this proposition are similar to those in~\cite{Hor76}. The first is obvious. For the second, we refer to~\cite{Hor76}. Properties $3.$ and $4.$ are quite trivial. 
We prove the last property. By~\eqref{Holdernorm},   
\begin{eqnarray}
\label{proofHoldercomp}
|f \circ z|_{C^\sigma} &\le& |f|_{C^0} + |(\nabla f\circ z)^T \nabla z |_{C^{\sigma -1}},
\end{eqnarray}
where $T$ stands for the transpose and $(\nabla f\circ z)^T \nabla z$ is the vector having $i$ component equal to $\left((\nabla f\circ z)^T \nabla z\right)_i = \nabla f\circ z \cdot \partial_{x_i} z$. Thanks to the  property $2.$
\begin{eqnarray*}
|f \circ z|_{C^\sigma} &\le& |f|_{C^0} + |(\nabla f\circ z)^T \nabla z |_{C^{\sigma -1}}\\
&\le& |f|_{C^0} +C(\sigma)|\nabla f\circ z|_{C^{\sigma -1}}|\nabla z |_{C^0} + C(\sigma)|\nabla f\circ z|_{C^0}|\nabla z |_{C^{\sigma -1}}.
\end{eqnarray*} 

The last term of the latter is bounded by $|f|_{C^1}|\nabla z |_{C^{\sigma -1}}$, it remains to estimate $|\nabla f\circ z|_{C^{\sigma -1}}|\nabla z |_{C^0}$. If $\sigma \le 2$, $|\nabla f \circ z|_{C^{\sigma-1}} \le | f|_{C^\sigma}|\nabla z|^{\sigma-1}_{C^0}+ |f|_{C^1}$ thanks to $4.$. Then 
\begin{eqnarray*}
|\nabla f \circ z|_{C^{\sigma-1}} |\nabla z|_{C^0}&\le& C(\sigma) \left(| f|_{C^\sigma}|\nabla z|^\sigma_{C^0}+ |f|_{C^1}|\nabla z|_{C^0}\right)\\
&\le& C(\sigma) \left(| f|_{C^\sigma}|\nabla z|^\sigma_{C^0}+ |f|_{C^1}|\nabla z|_{C^{\sigma -1}}\right),
\end{eqnarray*}
whence the property holds in this case. If $\sigma >2$, assuming that $5.$ is already proven for $\sigma -1$, we find
\begin{equation*}
|\nabla f \circ z|_{C^{\sigma-1}} |\nabla z|_{C^0}\le C(\sigma)\left(|\nabla f|_{C^{\sigma-1}}|\nabla z|^\sigma_{C^0}+ |f|_{C^2}|\nabla z|_{C^{\sigma-2}}|\nabla z|_{C^0} + |f|_{C^1}|\nabla z|_{C^0}\right).
\end{equation*}
It remains to find a good estimate for the central term. By Proposition \ref{convexity}
\begin{eqnarray*}
|f|_{C^2}|\nabla z|_{C^{\sigma-2}}|\nabla z|_{C^0} &\le& C(\sigma) \left(|f|^{\sigma-2 \over \sigma-1}_{C^1}|f|_{C^\sigma}^{1\over \sigma -1}\right) \left(|\nabla z|^{1\over \sigma -1}_{C^0}|\nabla z|_{C^{\sigma-1}}^{\sigma-2\over \sigma -1}\right)|\nabla z|_{C^0} \\
&\le&C(\sigma)  \left(|f|_{C^1}|\nabla z|_{C^{\sigma-1}}\right)^{\sigma-2 \over \sigma-1}  \left(|f|_{C^\sigma}|\nabla z|_{C^0}^\sigma\right)^{1 \over \sigma-1},
\end{eqnarray*}
since $a^\lambda b^{1-\lambda} \le C(a+b)$ for $0 < \lambda < 1$, we have the claim. 
\end{proof}

\section{Real analytic classes of functions}\label{B}

This section will collect some well-known facts about real analytic functions. For some $s > 0$, we begin with the introduction of complex domains 
\begin{eqnarray*}
\T_s^n &:=& \{q \in \C^n/\Z^n : | \operatorname{Im}(q)| \le s \} \\
B_s &:=&\{p \in \C^n : |p| \le s\},
\end{eqnarray*}
with $\T^n = \R^n/\Z^n$ and $B \subset \R^n$ a sufficiently large neighborhood of the origin. Let $D$ be equal to $\T^n \times B$ or $\T^n$ and we consider a real analytic function in a neighborhood of $D$ 
\begin{equation*}
f: D \to \R.
\end{equation*}
Let $D_s$ be equal to $\T_s^n \times B_s$ or $\T_s^n$, for a suitable small $s$. It is known that $f$ extends to a function 
\begin{equation*}
f: D_s \to \C
\end{equation*}
that is real, holomorphic and bounded. We define the following norm
\begin{equation*}
|f|_s = \sup_{z \in D_s}|f(z)|.
\end{equation*}
In the case of vector-valued functions $f = (f_1,...,f_n)$ with values in $\C^n$, we set $|f|_s = \max_i|f_i|_s$. Moreover, if $C=\{C_{ij}\}_{1 \le i,j\le n}$ is a $n \times n$ matrix, we let $|C|_s = \max_{ij}|C_{ij}|_s$. We define $\mathcal{A}_s$ as the space of such functions. The rest of this section is devoted to a series of general well-known properties.
\begin{proposition}
Let $f$, $g \in \mathcal{A}_s$, then the product $fg \in \mathcal{A}_s$ and 
\begin{equation*}
|fg|_s \le |f|_s|g|_s.
\end{equation*}
Let $f \in \mathcal{A}_s$ and $0 \le \sigma \le s$. Then $\partial_x f \in \mathcal{A}_s$ and we have
\begin{equation*}
|\partial_x f|_{s-\sigma} \le {1 \over \sigma}|f|_s.
\end{equation*}
Let $f \in \mathcal{A}_s$, $0 \le \sigma \le s$ and $\phi \in \mathcal{A}_{s-\sigma}$ such that $\phi: D_{s -\sigma} \to D_s$. Then $f \circ \phi \in \mathcal{A}_{s - \sigma}$ and 
\begin{equation*}
|f \circ \phi|_{s-\sigma}\le |f|_s.
\end{equation*}
\end{proposition}

\section{Banach spaces}\label{BanSp}
Here, we prove that the normed spaces introduced in Section \ref{PrelSettGD} are Banach spaces. Similarly, we have the claim for those defined in Section \ref{PSAGD}. 

Given $\upsilon \ge 0$, let $\mathbf{b}$ be a positive, decreasing, integrable function on $J_\upsilon$. We recall that
\begin{equation*}
\mathbf{\bar b}(t) = \int_t^{+\infty} \mathbf{b}(\tau) d\tau.
\end{equation*}
We may assume 
\begin{equation*}
\mathbf{b}(t) \le 1, \quad \mathbf{\bar b}(t) \le 1
\end{equation*}
for all $t \in J_\upsilon$. For fixed $\sigma \ge 1$ and an integer $k \ge 0$, we consider the following spaces $(\mathcal{G}, |\cdot|)$ and $(\mathcal{W}, \left\|\cdot \right\|)$  such that
\begin{align*}
&\mathcal{G} = \{g : \T^n \times J_\upsilon \to \R : g \in \mathcal{\bar S}^\upsilon_{\sigma, k} \hspace{2mm} \mbox{and}   \hspace{2mm} |g| =|g|^\upsilon_{\sigma + k, \mathbf{b}}< \infty\}\\
&\mathcal{W}^\sigma =\Big\{w : \T^n \times J_\upsilon  \to \R^n : w, \left(\nabla w \right)\Omega \in \mathcal{S}^\upsilon_{\sigma} \\
& \mbox{and}   \hspace{2mm}\left \|w\right \|_\sigma =\max \{|w|^\upsilon_{\sigma,\mathbf{\bar b}},|\left(\nabla w \right)\Omega|^\upsilon_{\sigma, \mathbf{b}}\} < \infty \Big\}
\end{align*}
We recall that, for all $(q, t)\in \T^n \times J_\upsilon$,  $\nabla w(q,t)\Omega = \partial_q w(q,t) \omega + \partial_t w(q,t)$ with $\omega \in \R^n$. Moreover, the spaces $\mathcal{S}^\upsilon_\sigma$ and $\mathcal{\bar S}^\upsilon_{\sigma, k}$ are defined respectively by Definition \ref{SGD} and Definition \ref{barSGD}, and the norm $|\cdot|^\upsilon_{\sigma, \mathbf{b}}$ is introduced by~\eqref{defnormGD}. 

We prove that these spaces are complete. We begin with the first. Let $\{g_d\}_{d \ge 0} \subset \mathcal{G}$ be a Cauchy sequence. This means that, for all $\varepsilon >0$ there exists $D \in \N$ such that for all $d$, $m \ge D$, $|g_d - g_m|^\upsilon_{\sigma +k , \mathbf{b}} < \varepsilon$. For all fixed $t \in J_\upsilon$, we claim that the sequence $\{g_d^t\}_{d \ge 0}$ contained in the Banach space $\left(C^{\sigma+k}(\T^n), |\cdot|_{C^{\sigma+k}}\right)$ is a Cauchy sequence. This is because for all fixed $t \in J_\upsilon$
\begin{equation*}
|g^t_d - g^t_m|_{C^{\sigma+k}} \le {|g^t_d - g^t_m|_{C^{\sigma+k}} \over \mathbf{b}(t)}  \le |g_d - g_m|^\upsilon_{\sigma +k , \mathbf{b}}.
\end{equation*}

Then, for all fixed $t \in J_\upsilon$, there exists $g^t \in C^{\sigma + k}$ such that 
\begin{equation*}
\lim_{d \to +\infty} |g^t_d - g^t|_{C^{\sigma +k}} = 0. 
\end{equation*}
We have to verify that $g \in \mathcal{G}$ (that is $g \in \mathcal{\bar S}^\upsilon_{\sigma, k}$ and $|g|^\upsilon_{\sigma +k , \mathbf{b}} <\infty$) and \newline$\lim_{d \to +\infty} |g_d - g|^\upsilon_{\sigma +k , \mathbf{b}} = 0$.

\textit{We prove that $g \in \mathcal{\bar S}^\upsilon_{\sigma, k}$}.  Obviously, for all fixed $t \in J_\upsilon$, $g^t \in C^{\sigma+k}(\T^n)$. It remains to verify that $\partial_q^i g \in C(\T^n \times J_\upsilon)$ for all $0 \le i\le k$. For all $(q_1, t_1)$, $(q_2, t_2) \in \T^n \times J_\upsilon$ 
\begin{eqnarray*}
|\partial_q^i g^{t_1}(q_1) - \partial_q^i g^{t_2}(q_2)| &\le& |\partial_q^i g_d^{t_1}(q_1) - \partial_q^i g^{t_1}(q_1)| + |\partial_q^i g_d^{t_1}(q_1) - \partial_q^i g_d^{t_2}(q_2)|\\
&+& |\partial_q^i g_d^{t_2}(q_2) - \partial_q^i g^{t_2}(q_2)|
\end{eqnarray*}
for all $0\le i \le k$. Now, for all $\varepsilon >0$ there exists $D \in \N$ such that, for all $d \ge D$, the first and the last term on the right-hand side of the latter are smaller than ${\varepsilon \over 3}$. This is because, for all fixed $t \in J_\upsilon$, $g_d^t$ converges to $g^t$ in the norm $C^{\sigma +k}$. Concerning the second term, we know that $g_d \in \mathcal{\bar S}^\upsilon_{\sigma, k}$. Hence, by the definition of $ \mathcal{\bar S}^\upsilon_{\sigma, k}$, $\partial_q^i g_d \in C(\T^n \times J_\upsilon)$ for all $0 \le i \le k$. Then, there exists $\delta >0$ such that if $|(q_1, t_1) - (q_2, t_2)|<\delta$ also the second term on the right-hand side of the latter is smaller than ${\varepsilon \over 3}$. This proves the claim. 

\textit{We prove that  $\displaystyle\lim_{d \to +\infty} |g_d - g|^\upsilon_{\sigma +k , \mathbf{b}} = 0$}. Let $g_{k_d}$ be a subsequence of $g_d$ such that 
\begin{equation*}
|g_{k_{d+1}} - g_{k_d}|^\upsilon_{\sigma+k, \mathbf{b}} < \left({1 \over 2}\right)^d
\end{equation*}
for all $d \in \N$. We claim that it suffices to prove the above property for $g_{k_d}$. Indeed, we assume that  $\displaystyle\lim_{d \to +\infty} |g_{k_d} - g|^\upsilon_{\sigma +k , \mathbf{b}} = 0$. Then, for all $d \in \N$
\begin{equation*}
|g_d - g|^\upsilon_{\sigma+k, \mathbf{b}} \le |g_d - g_{k_d}|^\upsilon_{\sigma+k, \mathbf{b}}  + |g_{k_d} - g|^\upsilon_{\sigma+k, \mathbf{b}}.
\end{equation*}
Therefore, for all $\varepsilon >0$, there exists $D \in \N$ such that, $|g_d - g_{k_d}|^\upsilon_{\sigma+k, \mathbf{b}} <{\varepsilon \over 2}$ and $|g_{k_d} - g|^\upsilon_{\sigma+k, \mathbf{b}}<{\varepsilon \over 2}$ for all $d \ge D$. Because $g_d$ is a Cauchy sequence, we have the first inequality. The second follows because we assumed that $g_{k_d}$ converges to $g$ in the norm $|\cdot|^\upsilon_{\sigma +k , \mathbf{b}}$. This implies $\displaystyle\lim_{d \to +\infty} |g_d - g|^\upsilon_{\sigma +k , \mathbf{b}} = 0$ and hence the claim.

Now, for all fixed $t \in J_\upsilon$
\begin{eqnarray*}
{|g^t_{k_d} - g^t|_{C^{\sigma+k}} \over \mathbf{b}(t)} \le \sum_{i=d}^{+\infty}{|g^t_{k_i} - g_{k_{i+1}}^t|_{C^{\sigma+k}} \over \mathbf{b}(t)} \le \sum_{i=d}^{+\infty}|g_{k_i} - g_{k_{i+1}}|^\upsilon_{\sigma + k,\mathbf{b}} \le  \sum_{i=d}^{+\infty} \left({1 \over 2}\right)^i = 2\left({1 \over 2}\right)^d
\end{eqnarray*}
and hence, taking the sup for all $t \in J_\upsilon$, we obtain
\begin{equation*}
|g_{k_d} - g|^\upsilon_{\sigma+k, \mathbf{b}} \le  2\left({1 \over 2}\right)^d.
\end{equation*}
Then, for every $\varepsilon >0$ there exists $D \in \N$ such that $|g_{k_d} - g|^\upsilon_{\sigma+k, \mathbf{b}} < \varepsilon$ for all $d \ge D$. 

\textit{We prove that  $ |g|^\upsilon_{\sigma +k , \mathbf{b}} < \infty$}. For all $d \in \N$, we can estimate $|g|^\upsilon_{\sigma +k , \mathbf{b}}$ as follows
\begin{equation*}
|g|^\upsilon_{\sigma +k , \mathbf{b}} \le |g_d - g|^\upsilon_{\sigma +k , \mathbf{b}} + |g_d|^\upsilon_{\sigma +k , \mathbf{b}}.
\end{equation*}
For $d$ sufficiently large,  $|g_d - g|^\upsilon_{\sigma +k , \mathbf{b}} <\infty$ because $\displaystyle\lim_{d \to +\infty} |g_d - g|^\upsilon_{\sigma +k , \mathbf{b}} = 0$. Moreover, $|g_d|^\upsilon_{\sigma +k , \mathbf{b}}<\infty$ because $g_d \in \mathcal{G}$. Then $(\mathcal{G}, |\cdot|)$ is a Banach space. 

In the second part of this section we prove that $(\mathcal{W}, \left\|\cdot \right\|)$ is a Banach space. 
Let $\{w_d\}_{d \ge 0} \subset \mathcal{W}$ be a Cauchy sequence. Similarly to the previous case, there exist $w\in \mathcal{S}^\upsilon_{\sigma}$ and $f \in \mathcal{S}^\upsilon_\sigma$ such that 
\begin{equation}
\label{limitv}
\lim_{d \to +\infty}|w_d - w|^\upsilon_{\sigma, \mathbf{\bar b}} = 0, \quad \lim_{d \to +\infty}|\left(\nabla w_d \right) \Omega- f|^\upsilon_{\sigma, \mathbf{b}} = 0.
\end{equation}
We have to verify that $\nabla w(q,t) \Omega = f(q,t)$ for all $(q,t) \in \T^n \times J_\upsilon$. Let us denote  $z = (q,t)$ and we recall that $\Omega = (\omega,1)$. We will prove that, for all $\varepsilon > 0$, there exists $\delta > 0$ such that 
\begin{equation*}
\left|{w(z + \tau\Omega) - w(z) \over \tau} - f(z)  \right| < \varepsilon
\end{equation*}
for all $|\tau| < \delta$. Thanks to the triangle inequality
\begin{eqnarray*}
\left|{w(z + \tau \Omega) - w(z) \over \tau} - f(z)  \right| &<& \left|{w_d(z + \tau\Omega) - w_d(z) \over \tau} - {w(z + \tau  \Omega) - w(z) \over \tau}  \right|\\
&+& \left|{w_d(z + \tau\Omega) - w_d(z) \over \tau} -  \nabla w_d(z)   \Omega\right|\\
&+& \left| \nabla w_d(z) \Omega - f(z)\right|.
\end{eqnarray*}
By~\eqref{limitv} there exists $D>0$, depending on $\varepsilon$ and $\tau$, such that the first and the third terms on the right-hand side of the latter are smaller than ${\varepsilon \over 3}$ for all $d \ge D$. Now, thanks to Taylor's formula, we can rewrite the second term on the right-hand side of the latter as follows
\begin{equation*}
\left|{w_d(z + \tau\Omega) - w_d(z) \over \tau} -  \nabla w_d(z)   \Omega\right| = \left|\int_0^1 \nabla w_d(z +\tau s \Omega)  \Omega -  \nabla w_d(z) \Omega ds  \right|
\end{equation*}
and using  the triangle inequality
\begin{eqnarray*}
\left|\int_0^1 \nabla w_d(z +\tau s \Omega)  \Omega -  \nabla w_d(z) \Omega ds  \right| &\le& \int_0^1 \left|  \nabla w_d(z +\tau s \Omega)  \Omega- f(z +\tau s \Omega)\right|d s \\
&+& \int_0^1 \left|f(z +\tau s \Omega) - f(z) \right| ds \\
&+& \int_0^1 \left|f(z) - \nabla w_d(z) \Omega\right| ds.
\end{eqnarray*}
We know that $f$ is continuous, then there exists $\delta$ such that for all $|\tau| < \delta$ the second term on the right-hand side of the latter is smaller than ${\varepsilon \over 9}$. Since the uniform convergence of $ \left(\nabla w_d\right) \Omega$ there exists $D>0$, depending on $\varepsilon$ and $\tau$, such that the first and the third terms on the right-hand side of the latter are smaller than ${\varepsilon \over 9}$.  This concludes the proof. 

\section*{Acknowledgement}

\textit{This project has received funding from the European Union’s Horizon 2020 research and innovation programme under the Marie Skłodowska-Curie grant agreement No 754362}.

\bibliographystyle{amsalpha}
\bibliography{Art1}

\providecommand{\bysame}{\leavevmode\hbox to3em{\hrulefill}\thinspace}
\providecommand{\MR}{\relax\ifhmode\unskip\space\fi MR }
\providecommand{\MRhref}[2]{%
  \href{http://www.ams.org/mathscinet-getitem?mr=#1}{#2}
}
\providecommand{\href}[2]{#2}
\begin{thebibliography}{CdlL15}

\bibitem[Arn63]{Arn63a}
V.~I. Arnold, \emph{Proof of a theorem of {A}. {N}. {K}olmogorov on the
  preservation of conditionally periodic motions under a small perturbation of
  the {H}amiltonian}, Uspehi Mat. Nauk \textbf{18} (1963), no.~5 (113), 13--40.
  \MR{0163025}

\bibitem[CdlL15]{CdlL15}
Marta Canadell and Rafael de~la Llave, \emph{K{AM} tori and whiskered invariant
  tori for non-autonomous systems}, Phys. D \textbf{310} (2015), 104--113.
  \MR{3396853}

\bibitem[FW14]{FW14}
Alessandro Fortunati and Stephen Wiggins, \emph{Persistence of {D}iophantine
  flows for quadratic nearly integrable {H}amiltonians under slowly decaying
  aperiodic time dependence}, Regul. Chaotic Dyn. \textbf{19} (2014), no.~5,
  586--600. \MR{3266829}

\bibitem[Hö76]{Hor76}
Lars Hörmander, \emph{The boundary problems of physical geodesy}, Arch.
  Rational Mech. Anal. \textbf{62} (1976), no.~1, 1--52. \MR{602181}

\bibitem[Kol54]{Kol54}
A.~N. Kolmogorov, \emph{On conservation of conditionally periodic motions for a
  small change in {H}amilton's function}, Dokl. Akad. Nauk SSSR (N.S.)
  \textbf{98} (1954), 527--530. \MR{0068687}

\bibitem[Mos62]{Mos62}
J.~Moser, \emph{On invariant curves of area-preserving mappings of an annulus},
  Nachr. Akad. Wiss. G\"{o}ttingen Math.-Phys. Kl. II \textbf{1962} (1962),
  1--20. \MR{147741}

\bibitem[Sca22a]{Sca22c}
D.~Scarcella, \emph{Biaymptotically quasiperiodic solutions for time-dependent
  hamiltonians}, In preparation (2022).

\bibitem[Sca22b]{Sca22b}
\bysame, \emph{Weakly asymptotically quasiperiodic solutions for time-dependent
  hamiltonians with a view to celestial mechanics}, In preparation (2022).

\end{thebibliography}

\end{document}